\documentclass[aop,MSNbibl,nameyear,dvips]{arximspdf}
\usepackage{mathrsfs}
\usepackage{mathbh}

\doi{10.1214/13-AOP895} 
\volume{42}
\issue{5}
\pubyear{2014}
\firstpage{1911}
\lastpage{1951}

\makeatletter
\newcommand{\mvert}{\vert}
\renewcommand{\mid}{\vert}

\newcommand{\rrVert}{\Vert}
\newcommand{\llVert}{\Vert}

\newproclaim{definition}{Definition}
\newproclaim{remark}[definition]{Remark}
\newproclaim{example}[definition]{Example}
\newtheorem{theorem}[definition]{Theorem}
\newtheorem{lemma}[definition]{Lemma}
\newtheorem{proposition}[definition]{Proposition}
\newtheorem{corollary}[definition]{Corollary}

\newcommand{\dist}{\operatorname{dist}}
\newcommand{\Exp}{\operatorname{Exp}}
\newcommand{\gibbs}{\mathrm{Gibbs}}
\newcommand{\Po}{\mathrm{Po}}
\newcommand{\tXXi}{\widetilde{\Xi}}

\newcommand{\Leb}{\mathrm{Leb}}
\newcommand{\SSS}{(\mathrm{S})}
\newcommand{\RS}{(\mathrm{RS})}
\newcommand{\UB}{(\mathrm{UB})}
\newcommand{\RC}{(\mathrm{RC})}
\newcommand{\IR}{(\mathrm{IR})}
\newcommand{\pip}{\mathrm{PIP}}
\makeatother

\begin{document}
\begin{frontmatter}

\title{Gibbs point process approximation: Total variation bounds using
Stein's method\thanksref{T1}}
\runtitle{Gibbs point process approximation}

\begin{aug}
\author{\fnms{Dominic} \snm{Schuhmacher}\corref{}\ead[label=e1]{dominic.schuhmacher@mathematik.uni-goettingen.de}}
\and
\author{\fnms{Kaspar} \snm{Stucki}\ead[label=e2]{kaspar.stucki@mathematik.uni-goettingen.de}}
\runauthor{D. Schuhmacher and K. Stucki}
\affiliation{University of Bern and University of G\"ottingen}
\address{Institute for Mathematical Stochastics\\
Georg-August-Universit\"at G\"ottingen\\
Goldschmidtstra{\ss}e 7\\
37077 G\"ottingen\\
Germany\\
\printead{e1}\\
\phantom{E-mail:\ }\printead*{e2}} 
\end{aug}
\thankstext{T1}{Supported by Swiss National Science Foundation Grant
200021-137527.}

\received{\smonth{7} \syear{2012}}
\revised{\smonth{4} \syear{2013}}

%
\begin{abstract}
We obtain upper bounds for the total variation distance between the
distributions of two Gibbs point processes in a very general setting.
Applications are provided to various well-known processes and settings
from spatial statistics and statistical physics, including the
comparison of two Lennard--Jones processes, hard core approximation of
an area interaction process and the approximation of lattice processes
by a continuous Gibbs process.

Our proof of the main results is based on Stein's method. We construct
an explicit coupling between two spatial birth--death processes to
obtain Stein factors, and employ the Georgii--Nguyen--Zessin equation
for the total bound.
\end{abstract}

%
\begin{keyword}[class=AMS]
\kwd[Primary ]{60G55}
\kwd[; secondary ]{60J75}
\kwd{82B21}
\end{keyword}
\begin{keyword}
\kwd{Conditional intensity}
\kwd{pairwise interaction process}
\kwd{birth--death process}
\kwd{Stein's method}
\kwd{total variation distance}
\end{keyword}

\end{frontmatter}

\section{Introduction}\label{secintroduction}

Gibbs processes form one of the most important classes of point
processes in spatial statistics that may incorporate dependence between
the points [\citet{moellerwaage04}, Chapter~6]. They are
furthermore, mainly in the special guise of pairwise interaction
processes, one of the building blocks of modern statistical physics
[\citet{ruelle69}].

Up to the somewhat technical condition of hereditarity (see Section~\ref{sectechpre}),
a Gibbs process on a compact metric space $\mathcal{X}$ is simply a
point process whose distribution is absolutely continuous with respect
to a ``standard'' Poisson process distribution. It is thus a natural
counterpart in the point process world to a real-valued random variable
that has a density with respect to some natural reference measure.
A~notorious difficulty with Gibbs processes is that in most cases of
interest their densities can only be specified up to normalizing
constants, which typically renders explicit calculations, for example,
of the total variation distance between two such processes, difficult.

In the current paper we give for the first time a comprehensive theorem
about upper bounds on the total variation distance between Gibbs
process distributions in a very general setting. These bounds provide
natural rates of convergence in many asymptotic settings, and include
explicit constants, which are small if one of the Gibbs processes is
not too far away from a Poisson process.

For the important special case of bounding the distance between two
pairwise interaction processes $\Xi_1$ and $\Xi_2$ on $\mathcal
{X}\subset\mathbb{R}^D$ with densities proportional to $\beta
^{\vert \xi\vert} \prod_{\{x,y\} \subset\xi} \varphi_1(x-y)$ and
$\beta
^{\vert \xi\vert} \prod_{\{x,y\} \subset\xi} \varphi_2(x-y)$,
respectively, where $\varphi_1$ and $\varphi_2$ are pairwise
interaction functions that are bounded by one (inhibitory case), a
consequence of our results is that there is an explicitly computable
constant $C=C(\beta,\varphi_2)>0$ such that
%
%
\begin{equation}
\label{eqintropip} d_{\mathrm{TV}} \bigl( \mathscr{L}(\Xi_1),
\mathscr{L}(\Xi_2) \bigr) \leq C \llVert\varphi_1-
\varphi_2 \rrVert_{L^1}.
\end{equation}
If we relax the condition that the pairwise interaction functions are
bounded by one and require suitable stability conditions for $\Xi_1$
and $\Xi_2$ instead, we still obtain
%
%
\begin{equation}
d_{\mathrm{TV}} \bigl( \mathscr{L}(\Xi_1), \mathscr{L}(
\Xi_2) \bigr) \leq C(\varepsilon) \llVert\varphi_1-
\varphi_2 \rrVert_{L^1} + \varepsilon,
\end{equation}
where $\varepsilon$ can be chosen arbitrarily small, causing a bigger
$C(\varepsilon)$.
We give more explicit examples for Strauss, bi-scale Strauss, and
Lennard--Jones type processes in Sections~\ref{secdist} and~\ref
{secapplications}.

For our proof of the main results we develop Stein's method for Gibbs
process approximation. Using the generator approach by~\citet
{barbour88}, we re-express the total variation distance in terms of the
infinitesimal generator of a spatial birth--death process (SBDP) whose
stationary distribution is one of the Gibbs process distributions
involved. An upper bound is then obtained by constructing an explicit
coupling of such SBDPs in order to obtain the so-called Stein factor
and applying the Georgii--Nguyen--Zessin equation.

Previously Stein's method has been applied very successfully for
Poisson process approximation; see \citet{bb92}, \citet
{chenxia04} and \citet{schumi4}. Other notable developments in the
domain of point process approximation concentrate on compound Poisson
process approximation [\citet{barman02}] and on approximation by
certain point processes whose points are i.i.d. given their total
number [called polynomial birth--death proceses by the authors; see
\citet{xiazhang12}]. In the latter article the authors give
substantially improved bounds when replacing approximating Poisson or
Compound Poisson processes by their new processes. However, these new
processes are by no means flexible enough to approximate processes
typically encountered in spatial statistics, where truely local point
interactions take place, such as mutual inhibition up to a certain
(nonnegligible) distance.

In \citet{steinintro} the editors write in the preface:
``Point process approximation, other than in the Poisson context, is largely
unexplored.'' This statement still remains mostly true today, and the
present paper makes a substantial contribution in order to change this.

Apart from approaches by Stein's method, the authors are not aware of
any publications that give bounds for a probability metric between
Gibbs processes in any generality. There is, however, related work by
\citet{moeller89} and \citet{daipraposta12}, where
convergence rates
for distances between an SBDP and its stationary point process
distribution were considered.

The plan of the paper is as follows. We start out in Section~\ref
{sectechpre} by giving the necessary definitions and notation,
including a somewhat longer introduction to Gibbs and pairwise
interaction processes. Section~\ref{secdist} contains the main
results. While Section~\ref{ssecmain} treats the common case where
the approximating Gibbs process satisfies a stronger stability
condition, Sections~\ref{ssecnotS} and~\ref{ssecg-pip} lay out a
strategy and give concrete results under very general conditions.
Simpler examples are scattered throughout Section~\ref{secdist},
while Section~\ref{secapplications} looks at the three more involved
applications mentioned in the abstract. In Section~\ref{seccoupling}
we discuss spatial birth--death processes and present the coupling
needed for obtaining the Stein factors, and in Section~\ref
{secSteinforGibbs} we develop Stein's method for Gibbs process
approximation and give the proofs of the main results. The paper
finishes by an \hyperref[app]{Appendix} that justifies the reduction of our main proofs
to a state space with diffuse reference measure $\bolds{\alpha}$.

\section{Prerequisites}
\label{sectechpre}

Let $(\mathcal{X}, d)$ be a compact metric space, which serves
as the state space for all our point processes. We equip $\mathcal{X}$
with its Borel $\sigma$-algebra $\mathcal{B}=\mathcal{B}(\mathcal
{X})$. Let $\bolds{\alpha}\neq0$ be a fixed finite reference measure on
$(\mathcal{X},\mathcal{B})$. If $\mathcal{X}$ has a suitable group
structure, $\bolds{\alpha}$ is typically chosen to be the Haar
measure. If
$\mathcal{X}\subset\mathbb{R}^D$, we tacitly use Lebesgue measure
and write $\vert A \vert:= \Leb^D(A)$. Also, unless specified otherwise,
we assume $d$ to be the Euclidean metric in this case and write
$\alpha_D = \pi^{D/2} / \Gamma(D/2+1)$ for the volume of the
unit ball. 

Denote by $(\mathfrak{N},\mathcal{N})$ the space of finite counting
measures (``point configurations'') on $\mathcal{X}$ equipped with its
canonical $\sigma$-algebra; see \citet{kallenberg86}, Section~1.1.
For any $\xi\in\mathfrak{N}$ write $\vert \xi\vert = \xi
(\mathcal
{X})$ for its total number of points. A~\emph{point process} is simply
a random element of $\mathfrak{N}$.

For a finite measure $\bolds{\lambda}$ on $\mathcal{X}$ recall that
a point
process $\Pi$ is called a \emph{Poisson process} with \emph
{intensity} measure $\bolds{\lambda}$ if the point counts $\Pi
(A_i)$, $1
\leq i \leq n$, are independent $\Po(\bolds{\lambda}(A_i))$-distributed
random variables for any $n \geq1$ and any pairwise disjoint sets
$A_1,\ldots,A_n \in\mathcal{B}$. It is a well-known fact that such a
Poisson process may be constructed as $\Pi= \sum_{i=1}^N \delta
_{X_i}$, where $N$ is a $\Po(\bolds{\lambda}(\mathcal
{X}))$-distributed random variable, and $X_i$ are i.i.d. random
elements of $\mathcal{X}$ with distribution $\bolds{\lambda}(\cdot
)/\bolds{\lambda}(\mathcal{X})$ that are independent of $N$. We
denote the
Poisson process distribution with intensity measure $\bolds{\alpha}$ by
$\Po_1$. We will make extensive use of the fact that\vadjust{\goodbreak} for $\Pi
\sim \Po_1$ and any measurable function $h \colon\mathfrak
{N}\to\mathbb{R}_{+}$ we have
%
%
\begin{eqnarray}
\label{eqpkf} 
\mathbb{E}h(\Pi) &=& \int_{\mathfrak{N}} h(
\xi) \Po_1(d\xi)
\nonumber\\[-8pt]\\[-8pt]
&=& e^{-\bolds{\alpha}(\mathcal{X})} \sum_{k=0}^{\infty}
\frac
{1}{k!} \int_\mathcal{X}\cdots\int_\mathcal{X}h
\Biggl(\sum_{i=1}^k \delta
_{x_i} \Biggr) \bolds{\alpha}(dx_1)\cdots\bolds{
\alpha}(dx_k),\nonumber
\end{eqnarray}
where we interpret the summand for $k=0$ as $h(\varnothing)$, writing
$\varnothing$ for the empty point configuration. Equation~(\ref
{eqpkf}) is obtained by conditioning on the total number of points of
$\Pi$. Note that a similar version may also be found in \citet{moellerwaage04}, Proposition~3.1(ii).

\subsection{Gibbs processes}

We give the definition of a Gibbs point process from spatial
statistics. Use the natural partial order on $\mathfrak{N}$, by which
$\xi\leq\eta$ if and only if $\xi(A) \leq\eta(A)$ for every $A
\in\mathcal{B}$. We call a function $u \colon\mathfrak{N}\to
\mathbb{R}_{+}$ \emph{hereditary} if for any $\xi, \eta\in
\mathfrak{N}$ with $\xi\leq\eta$, we have that $u(\xi) = 0$
implies $u(\eta) = 0$.
%
%
\begin{definition}
A point process $\Xi$ on $\mathcal{X}$ is called a \emph{Gibbs
process} if it has a~hereditary density $u$ with respect to $\Po_1$.
\end{definition}
Gibbs processes form important models in spatial statistics. The
specification of a density allows us to model point interactions in a
simple and intuitive way. Under rather flexible conditions, for
example, log-linearity in the parameters, formal inference (partly
based on numerical methods) is possible. See \citet{moellerwaage04}, Chapter~9.

It will be convenient to identify a Gibbs process by its conditional intensity.
%
%
\begin{definition}
Let $\Xi$ be a Gibbs process with density $u$. We call the function
$\lambda(\cdot\mvert\cdot) \colon\mathcal{X}\times\mathfrak
{N}\to\mathbb{R}_{+}$,
%
%
\begin{equation}
\label{eqcif} \lambda(x \mvert\xi) = \frac{u(\xi+ \delta
_{x})}{u(\xi)},
\end{equation}
the \emph{conditional intensity (function)} of $\Xi$. For this
definition we use the convention that $0/0 = 0$.
\end{definition}

Note that other definitions of the conditional intensity in the
literature may differ at pairs $(x,\xi)$ with $x \in\xi$. It is well
known and with the help of equation~(\ref{eqpkf}), straightforward to
check that the conditional intensity is the $\bolds{\alpha}\otimes
\mathscr
{L}(\Xi)$-almost everywhere unique product measurable function that
satisfies the \emph{Georgii--Nguyen--Zessin equation}
%
%
\begin{equation}
\label{eqgnz} \mathbb{E} \biggl( \int_{\mathcal{X}} h(x, \Xi-
\delta_x) \Xi(d x) \biggr) = \int_{\mathcal{X}}
\mathbb{E} \bigl( h(x, \Xi) \lambda(x \mvert\Xi) \bigr) \bolds
{\alpha}(d x)
\end{equation}
for every measurable $h \colon\mathcal{X}\times\mathfrak{N}\to
\mathbb{R}_{+}$.

A Gibbs process is usually specified via an unnormalized density
$\tilde{u}$ that is shown to be $\Po_1$-integrable. Typically
the integral and hence the normalized density $u$~cannot be computed
explicitly. On the other hand the conditional intensity can be
calculated simply as
%
%
\begin{equation}
\label{eqcifohne} \lambda(x \mvert\xi) = \frac{\tilde{u}(\xi+
\delta_{x})}{\tilde
{u}(\xi)}
\end{equation}
and has a nice intuitive interpretation as the infinitesimal
probability that $\Xi$ produces a (further) point around $x$ given it
produces at least the point configuration~$\xi$. Also it determines
the Gibbs process distribution completely, since an unnormalized
density $\tilde{u}$ can be recovered recursively for increasingly large
point configurations by employing~(\ref{eqcifohne}).
We denote by $\gibbs(\lambda)$ the distribution of the Gibbs process
with conditional intensity $\lambda= \lambda(\cdot\mvert\cdot)$.

The measure $\bolds{\lambda}$ given by $\bolds{\lambda}(A)={\mathbb
E}(\Xi(A))$ for any $A
\in\mathcal{B}$ is called the \emph{intensity measure} of $\Xi$,
provided that it is finite. By equation~(\ref{eqgnz}) we have
\[
\bolds{\lambda}(A)={\mathbb E} \biggl(\int_\mathcal{X}
\mathbh{1}\{ x\in A\} \Xi(dx) \biggr)=\int_\mathcal{X}\mathbh{1}
\{x\in A\} {\mathbb E}\bigl(\lambda(x\mid\Xi)\bigr) \bolds{\alpha}(dx),
\]
that is, $\bolds{\lambda}$ is absolutely continuous with respect to
$\bolds{\alpha}
$. We call its density $\lambda(x)={\mathbb E}(\lambda(x\mid\Xi))$ the
\emph{intensity (function)} of $\Xi$. We use notation $\lambda(\cdot
)$ and $\lambda(\cdot\mvert\cdot)$ to distinguish the intensity and
the conditional intensity if necessary.

In the main part of the paper we distinguish between the \emph
{approximated Gibbs process} $\Xi$ with a general conditional
intensity $\nu$, and the \emph{approximating Gibbs process} $\mathrm{H}$,
whose conditional intensity $\lambda$ will typically (except in
Sections~\mbox{\ref{ssecnotS}--\ref{sseclenny-jones}}) satisfy the
stability condition
{\renewcommand{\theequation}{S}{
\begin{equation}
\sup_{\xi\in\mathfrak{N}} \int_{\mathcal{X}}\lambda(x\mvert\xi)
\bolds{\alpha}(dx) < \infty.
\end{equation}}}\setcounter{equation}{6}%
Note that this condition follows from the \emph{local stability condition}
\[
\lambda(x \mvert\xi) \leq\psi^{*}(x)
\]
for an integrable function $\psi^{*} \colon\mathcal{X}\to\mathbb
{R}_{+}$. Local stability is satisfied for many point process
distributions traditionally used in spatial statistics. See \citet
{moellerwaage04}, page~84ff.

\subsection{Pairwise interaction processes}

A special type of Gibbs processes that are noteworthy both for their
relative simplicity and their abundant use in statistical physics are
the pairwise interaction processes. We treat distances between such
processes in detail in Sections~\ref{secdist} and~\ref{secapplications}.
%
%
\begin{definition}
A Gibbs process $\Xi$ on $\mathcal{X}$ is called a \emph{pairwise
interaction process} (PIP) if there exist $\beta:\mathcal{X}\to
\mathbb{R}_+$ and symmetric $\varphi:\mathcal{X}\times\mathcal
{X}\to\mathbb{R}_+$ such that $\Xi$ has the unnormalized density
\[
\tilde{u}(\xi)=\prod_{1 \leq i \leq n}\beta(x_i)
\prod_{1 \leq i <
j \leq
n}\varphi(x_i,x_j)
\]
for any $\xi= \sum_{i=1}^n \delta_{x_i} \in\mathfrak{N}$. The
normalizing constant is usually not analytically computable. We then
denote the distribution of $\Xi$ by $\pip(\beta,\varphi)$. The PIP
is called \emph{inhibitory} if $\varphi\leq1$. It is called \emph
{hard core} with radius $\delta> 0$ if $\varphi(x,y) = 0$ whenever
$d(x,y) \leq\delta$.
\end{definition}

The conditional intensity of $\Xi\sim\pip(\beta,\varphi)$ is
accordingly given by
\[
\lambda(x\mvert\xi)=\beta(x)\prod_{i = 1}^n
\varphi(x,x_i).
\]

For $\tilde{u}$ to be integrable with respect to $\Po_1$, it is by
equation~(\ref{eqpkf}) necessary that $\beta$~is integrable. For
inhibitory PIPs this is obviously also sufficient. The same holds for
hard core PIPs with bounded $\varphi$, because by the compactness
of~$\mathcal{X}$ their total number of points is almost surely
bounded. For more general PIPs the situation is not so simple; see
Example~\ref{exmsstrauss} for a special case. We will then assume the
following conditions:
\begin{longlist}[$\RS$]
\item[$\RS$] \textit{Ruelle stability}. There exist a constant $c^*$
and an integrable function $\psi^*$ such that $\tilde{u}(\xi)\le
c^*\prod_{i=1}^n \psi^*(x_i)$ for every $\xi= \sum_{i=1}^n \delta
_{x_i} \in
\mathfrak{N}$.
\end{longlist}
\begin{longlist}[$\UB$]
\item[$\UB$] \textit{Upper boundedness}. There exists a constant $C$
such that $\varphi(x,y) \le C$ for all $x,y\in\mathcal{X}$.
\end{longlist}
\begin{longlist}[$\RC$]
\item[$\RC$] \textit{Repulsion condition}. There exist $\delta>0$ and
$0\le\gamma\le1$ such that for all $x,y\in\mathcal{X}$ with
$d(x,y)\le\delta$ we have $\varphi(x,y) \le\gamma$.
\end{longlist}

Note that $\RS$ is the form of Ruelle stability commonly used in spatial
statistics; see \citet{moellerwaage04}. If we can choose $\psi^*(x)$
as $\beta(x)$ times a constant, we get the classical definition by
\citet{ruelle69}. In any case Ruelle stability ensures that the
unnormalized density $\tilde{u}$ is integrable.

If we can write the interaction function as $\varphi
(x,y)=e^{-V(x,y)}$, then $\UB$ is equivalent
to requiring that the potential $V$ is bounded from below, which is a
commonly used condition in statistical physics; see, for example, \citet{ruelle69}.

Furthermore, we introduce notation for the inner and outer ranges of
attractive interaction.
\begin{longlist}
\item[$\IR$] \textit{Interaction ranges}. Let $\delta\leq r<R$ be
constants such that for all $x,y\in\mathcal{X}$ with $d(x,y)\le
r$ or $d(x,y)> R$ we have $\varphi(x,y) \le1$.
\end{longlist}
Note that such constants always exist due to $\RC$ and the compactness
of~$\mathcal{X}$.

Strictly speaking only inhibitory PIPs satisfy $\SSS$. However, for our
purpose it is actually enough to require finiteness of both the
$\mathscr{L}(\Xi)$-$\operatorname{ess}\sup$ and the $\mathscr
{L}(\mathrm{H})$\mbox{-$\operatorname{ess}\sup$} instead of the supremum in $\SSS$. This
would also
admit comparisons
of arbitrary hard core PIPs. However, for the ease of presentation we
deal with hard core PIPs together with the more general PIPs in
Sections~\ref{ssecg-pip} and~\ref{sseclenny-jones}.

\subsection{Reduction to a diffuse reference measure $\bolds{\alpha}$}
\label{ssecdiffuse}

In the remainder of this paper we will tacitly assume that the
reference measure $\bolds{\alpha}$ is diffuse, that is, satisfies
$\bolds{\alpha}(\{
x\}) = 0$ for any $x \in\mathcal{X}$. This implies that the $\Po
_1$-process and the corresponding Gibbs processes are simple, that is,
with probability one do not have multiple points at a single location
in space. It is then convenient to interpret a point process as a
random finite set and use set notation, which is commonly done is
spatial statistics. Thus we may write, for example, $\xi\subset\eta$
instead of $\xi\leq\eta$, or in the density of a PIP $\prod_{x \in
\xi}\beta(x)\prod_{\{x,y\} \subset\xi}\varphi(x,y)$ instead of
$\prod_{1 \leq i \leq n}\beta(x_i)\prod_{1 \leq i < j \leq n}\varphi
(x_i,x_j)$.

In addition to simplifying the notation considerably by making points
identifiable by their location in space, assuming a diffuse $\bolds
{\alpha}$
also reduces the differences in various definitions of the conditional
intensity $\lambda(x \mvert\xi)$ for the case $x \in\xi$ to an
$\bolds{\alpha}\otimes \Po_1$-null set.

We show in the \hyperref[app]{Appendix} that our results in Section~\ref{secdist}
carry over to the nondiffuse case. Essentially this is seen by
extending the state space $\mathcal{X}$ to $\mathcal{X}\times[0,1]$
and considering $\bolds{\alpha}\otimes\Leb\vert_{[0,1]}$ as a new (always
diffuse) reference measure. This is based on an idea used in
\citet
{chenxia04}.

\section{Total variation bounds between Gibbs process distributions}
\label{secdist}

\subsection{Main results}
\label{ssecmain}

Define $\mathcal{F}_{\mathrm{TV}}$ as the set of all measurable
functions $f\colon\mathfrak{N}\to[0,1]$. Then, for two point processes
$\Xi$ and $\mathrm{H}$, the \emph{total variation distance} is
defined as
%
%
\begin{equation}
\label{eqdef-dtv} d_{\mathrm{TV}}\bigl(\mathscr{L}(\Xi),\mathscr
{L}(\mathrm{H})
\bigr)=\sup_{f\in\mathcal{F}_{\mathrm{TV}}}\bigl|{\mathbb E} f(\Xi
)-{\mathbb E}f(\mathrm{H})\bigr|.
\end{equation}
By a simple approximation argument this is equivalent to
%
%
\begin{equation}
d_{\mathrm{TV}}\bigl(\mathscr{L}(\Xi),\mathscr{L}(\mathrm
{H})\bigr)=\sup
_{A\in\mathcal
{N}}\bigl|{\mathbb P}(\Xi\in A)-{\mathbb P}(\mathrm{H}\in A)\bigr|.
\end{equation}
Denote by $\llVert \cdot\rrVert$ the total variation norm for
signed measures
on $\mathcal{X}$. Thus $\llVert \xi-\eta\rrVert$ for $\xi,\eta
\in
\mathfrak{N}$ is the total number of points appearing in one of the
point configurations, but not in the other.

Our main results are given as Theorems~\ref{thmmain2} and~\ref
{cordtv-pip}. The principal idea behind our proofs is a suitable
variant of Stein's method, which we develop in Section~\ref
{secSteinforGibbs}. The proofs themselves are deferred to Section~\ref
{secSteinforGibbs} as well.
%
%
\begin{theorem}
\label{thmmain2}
Let $\Xi\sim\gibbs(\nu)$ and $\mathrm{H}\sim\gibbs(\lambda)$ be Gibbs
processes. Suppose that $\mathrm{H}$ satisfies $\SSS$. Then there is a finite
constant $c_1(\lambda)$ such that
%
%
\begin{equation}
\label{eqthm-main} d_{\mathrm{TV}}\bigl(\mathscr{L}(\Xi),\mathscr
{L}(\mathrm{H})
\bigr) \le c_1(\lambda)\int_\mathcal{X}\mathbb{E}\bigl|
\nu(x\mvert\Xi)-\lambda(x \mvert\Xi)\bigr| \bolds{\alpha}(dx).
\end{equation}
More precisely, we have for any $n^*\in\mathbb{N}\cup\{\infty\}$ that
%
%
\begin{eqnarray}\label{eqc1}
c_1(\lambda) &\le& \bigl(n^*-1\bigr)! \biggl(
\frac{\varepsilon}{c} \biggr)^{n^*-1} \Biggl(\frac
{1}{c}\sum_{i=n^*}^\infty\frac{c^i}{i!}+\int
_0^c\frac{1}{s}\sum
_{i=n^*}^\infty\frac{s^i}{i!} \,ds \Biggr)
\nonumber\\[-8pt]\\[-8pt]
&&{} + \frac{1+\varepsilon}{\varepsilon}\sum_{i=1}^{n^*-1}
\frac{\varepsilon^{i}}{i},\nonumber
\end{eqnarray}
where
\[
\varepsilon = \sup_{\|\xi-\eta\|=1}\int_\mathcal{X}\bigl|
\lambda(x\mvert\xi)-\lambda(x \mvert\eta)\bigr| \bolds{\alpha}(dx) <
\infty
\]
and
\[
c=c\bigl(n^*\bigr)=\sup_{\|\xi-\eta\|\ge n^*}\int_\mathcal{X}\bigl|
\lambda(x\mvert\xi)-\lambda(x \mvert\eta)\bigr| \bolds{\alpha}(dx) <
\infty.
\]
If $n^{*}=\infty$, we interpret the long first summand in the upper
bound as $0$. For $\varepsilon=0$ and/or $c=0$ the upper bound is to be
understood in the limit sense.
\end{theorem}

%
\begin{remark}
The term $c_1(\lambda)$ has a special meaning in our proof and in the
theory of Stein's method in general; for its definition see
equation~(\ref{eqc1def}). It is usually referred to as the \emph
{(first) Stein factor}.
\end{remark}

%
\begin{remark}[{[Special cases for the bound on $c_1(\lambda)$]}]
\label{remspecialc1}
(a)~We often have a bound for $c$ which does not
depend on $n^*$. In
this case we choose $n^*=\lceil c/\varepsilon\rceil$, which turns out
to be optimal.

(b)~If $\varepsilon< 1$, we can choose $n^*=\infty$ and obtain
\[
c_1(\lambda) \le\frac{1+\varepsilon}{\varepsilon}\log\biggl
(\frac
{1}{1-\varepsilon}
\biggr) \le\frac{1+\varepsilon}{1-\varepsilon}.
\]
Conditions of the type $\varepsilon< 1$ are known in the statistical
physics literature as ``low activity, high temperature'' setting; see,
for example,~\citet{kl05}.

(c)~If $\mathrm{H}$ is a Poisson process, then $\lambda(x \mvert\xi) =
\lambda(x)$ does not depend on $\xi$. We then have $\varepsilon=0$, and
obtain $c_1(\lambda)=1$. Hence inequality~(\ref{eqc1}) contains the
bound on the first Stein factor given in Lemma~2.2(i) of \citet{bb92}
as a special case.
\end{remark}

%
\begin{remark}
\label{remubermain2}
The assumption that $\Xi$ is a Gibbs process is used only in the proof
of Theorem~\ref{thmmain2} for invoking the Georgii--Nguyen--Zessin
equation, and is of course not needed for bounding $c_1(\lambda)$. The
Georgii--Nguyen--Zessin equation may be generalized by replacing the
kernel $\nu(x \mvert\xi) \bolds{\alpha}(dx)$ with the
Papangelou kernel
$\bolds{\nu}(dx \mvert\xi)$ of $\Xi$ if the so-called
condition~($\Sigma
$) is satisfied, that is, if $\mathbb{P} (\Xi(A) = 0 \mvert \Xi
\vert_{A^c} ) > 0$ a.s. for every $A \in\mathcal{B}$. See
\citet{kallenberg86}, Section~13.2, or \citet{dvj08}, Section~15.6,
for details.

We may therefore generalize Theorem~\ref{thmmain2} as follows. Let
$\Xi$ be a point process that satisfies condition $(\Sigma)$ and has
Papangelou kernel $\bolds{\nu}$, and let $\mathrm{H}\sim\gibbs
(\lambda)$
satisfy condition $\SSS$. Then
\[
d_{\mathrm{TV}}\bigl(\mathscr{L}(\Xi),\mathscr{L}(\mathrm
{H})\bigr) \le
c_1(\lambda) \mathbb{E} \bigl\llVert\bolds{\nu}(dx\mvert\Xi)-
\lambda(x \mvert\Xi) \bolds{\alpha}(dx) \bigr\rrVert,
\]
where $c_1(\lambda)$ is as above, and \mbox{$\llVert \cdot\rrVert$} is
the total
variation norm for signed measures
on~$\mathcal{X}$.
\end{remark}

For inhibitory PIPs we obtain the following theorem, which relates the
total variation distance to the $L^1$-distance between the interaction
functions.

%
\begin{theorem}
\label{cordtv-pip}
Suppose that $\Xi\sim\pip(\beta,\varphi_1)$ and $\mathrm{H}\sim
\pip
(\beta,\varphi_2)$ are inhibitory. Let $\nu(y) = \mathbb{E}(\nu(y
\mvert
\Xi))$ denote the intensity of $\Xi$.
Then
%
%
\begin{eqnarray}\label{eqgeneralpipbound}
&& d_{\mathrm{TV}}\bigl(\mathscr{L}(\Xi),\mathscr{L}(\mathrm{H})\bigr)
\nonumber\\[-8pt]\\[-8pt]
&&\qquad \leq c_1(\lambda) \int_{\mathcal{X}} \int
_{\mathcal{X}} \beta(x) \nu(y) \bigl\vert\varphi_1(x,
y) - \varphi_2(x,y) \bigr\vert\bolds{\alpha} (d x) \bolds{
\alpha}(d y),\nonumber
\end{eqnarray}
where $c_1(\lambda)$ is bounded in inequality~(\ref{eqc1}) with
\[
c\leq\int_\mathcal{X}\beta(x) \bolds{\alpha}(dx) \quad
\mbox{and}\quad\varepsilon=\sup_{y\in\mathcal{X}}\int_\mathcal{X}
\beta(x) \bigl(1-\varphi_2(x,y)\bigr) \bolds{\alpha}(dx).
\]
In the case where $\mathcal{X}\subset\mathbb{R}^D$, $\beta$ is
constant, and $\varphi_i(x,y) = \varphi_i(x-y)$ depends only on the
difference, we obtain
%
%
\begin{equation}
\label{eqspecialpipbound} d_{\mathrm{TV}}\bigl(\mathscr{L}(\Xi
),\mathscr{L}(\mathrm{H})
\bigr) \leq c_1(\lambda) \beta\mathbb{E}\bigl(\vert\Xi\vert
\bigr) \int_{\mathbb{R}^D} \bigl\vert\varphi_1(x) -
\varphi_2(x) \bigr\vert \,d x.
\end{equation}
\end{theorem}
Note that $\nu(\cdot)$ can usually not be calculated explicitly, but
at least it can always be bounded by $\beta(\cdot)$. In particular,
inequality~(\ref{eqspecialpipbound}) implies inequality~(\ref{eqintropip}) in the \hyperref[secintroduction]{Introduction}. Note that better bounds on a
constant $\nu$ have been obtained in \citet{ss13pgfl}.

%
\begin{remark}
Our bounds on the total variation distance in Theorems \ref{thmmain2}~and~\ref{cordtv-pip} may be larger than one, in which case they give
no new information.
They are small if one of the processes is not too far away from a
Poisson process, and the conditional intensities (or the pairwise
interaction functions in the case of Theorem~\ref{cordtv-pip}) are
close in an $L^1$-sense. In what follows we are mainly interested in an
asymptotic setting, where, for example, the interaction function of a
$\pip$ converges to the interaction function of another $\pip$.
\end{remark}

If one of the processes is a Poisson process, we obtain both a slight
improvement and a very substantial generalization of the bounds
in~\citet{browngreig94}.
%
%
\begin{example}
Let $\Pi$ be a Poisson process with intensity function $\beta$, and
let $\Xi\sim\pip(\beta,\varphi)$ be inhibitory, denoting its
intensity function by $\nu$. By Theorem~\ref{cordtv-pip} and
Remark~\ref{remspecialc1}(c), we obtain
%
%
\begin{equation}
\label{eqpoissonpip} d_{\mathrm{TV}}\bigl(\mathscr{L}(\Xi
),\mathscr{L}(\Pi)\bigr)\le
\int_{\mathcal{X}} \int_{\mathcal{X}} \beta(x) \nu(y)
\bigl(1-\varphi(x,y)\bigr) \bolds{\alpha}(dx) \bolds{\alpha}(dy).
\end{equation}
The special case where $\mathcal{X}=[0,1]^D$ with torus convention,
$\bolds{\alpha}$ is Lebesgue measure,~and $\Xi$ is a stationary hard core
process with constant $\beta$ and $\varphi(x,y)=\mathbh{1}\{\llVert x-y
\rrVert >
r\}$ was considered in \citet{bb92} and \citet{browngreig94}, except
that these articles approximate by a Poisson process \mbox{$\widetilde{\Pi
}$} that has
the same intensity $\nu$ as $\Xi$. We obtain from (\ref
{eqpoissonpip}) that
%
%
\begin{equation}
\label{eqpoissonhc} d_{\mathrm{TV}}\bigl(\mathscr{L}(\Xi
),\mathscr{L}(\Pi)\bigr)
\leq\beta\nu\alpha_D r^D \leq\beta^2
\alpha_D r^D.
\end{equation}
The best bound in \citet{browngreig94}, namely inequality~(12), says
that under a somewhat complicated additional condition on the
parameters, we have
%
%
\begin{equation}
\label{eqbrowngreig} d_{\mathrm{TV}}\bigl(\mathscr{L}(\Xi
),\mathscr{L}(\widetilde{
\Pi})\bigr) \leq\biggl( 1 + \frac
{1}{2^D} \biggr) \nu^2
\alpha_D r^D.
\end{equation}
By a straightforward upper bound on $\beta$ [see \citet
{browngreig94}, inequality~(11)] our result (\ref{eqpoissonhc}) may
be bounded further to obtain
\[
d_{\mathrm{TV}}\bigl(\mathscr{L}(\Xi),\mathscr{L}(\Pi)\bigr)
\leq
\frac
{1}{1-\nu
\alpha_D r^D} \nu^2 \alpha_D r^D
\]
for $r < 1/(\nu\alpha_D)^{1/D}$, which holds without the additional
condition and is an asymptotic improvement over (\ref{eqbrowngreig})
by a factor of $2^D/(2^D+1)$ as $r \to0$.

This suggests that for small $r$ it is better to approximate $\Xi$ by
$\Pi$ than by $\widetilde{\Pi}$, both because we get a smaller bound
and because
the intensity of $\Pi$ is known explicitly from the parameters of $\Xi$.
\end{example}

%
\begin{example}
\label{exStrauss}
A PIP is called a \emph{Strauss process} if its interaction function
is given by
\[
\varphi(x,y)= %
\cases{ \gamma, &\quad if $d(x,y)\le R$,
\cr
1, &\quad
if $d(x,y)>R$} %
\]
for some constants $0\le\gamma\le1$ and $R>0$.
Let $\Xi$ and $\mathrm{H}$ be Strauss processes with constant $\beta
$ and
further parameters
$\gamma_1,R_1$ and $\gamma_2,R_2$, respectively, where $R_1>R_2$.
Denote by $\mathbb{B}(y,R)$ the closed ball in $\mathcal{X}$ with
center at
$y$ and radius $R$. Then by Theorem~\ref{cordtv-pip},
%
%
\begin{eqnarray}\label{eqdtv-Strauss}
\qquad && d_{\mathrm{TV}}\bigl(\mathscr{L}(\Xi),\mathscr{L}(\mathrm{H}) \bigr)\nonumber
\\
&&\qquad \le  c_1(\lambda){\mathbb E} \bigl(|\Xi|\bigr) \beta
\\
&&\quad\qquad {} \times \sup_{y\in\mathcal{X}} \bigl((1-\gamma_1)\bolds{\alpha}
\bigl(\mathbb{B} (y,R_1)\setminus\mathbb{B}(y,R_2)\bigr) +
|\gamma_1-\gamma_2|\bolds{\alpha}\bigl(
\mathbb{B}(y,R_2)\bigr) \bigr),\nonumber
\end{eqnarray}
where $c_1(\lambda)$ is bounded in inequality~(\ref{eqc1}) with
\[
\varepsilon=\beta(1-\gamma_2)\sup_{y\in\mathcal{X}}\bolds{
\alpha}\bigl(\mathbb{B} (y,R_2)\bigr) \quad\mbox{and}\quad c\leq
\beta
\bolds{\alpha}(\mathcal{X}).
\]
\end{example}

\subsection{Processes violating the stability condition $\SSS$}
\label{ssecnotS}
Many Gibbs processes satisfy condition $\SSS$, but there are some
important exceptions. In the present subsection we provide a technique
for treating these exceptions. In Section~\ref{ssecg-pip} we apply
this technique to general PIPs.

We call an event $A$ \emph{hereditary} if the corresponding
indicator function is hereditary, that is, if $\eta\in A$ implies $\xi
\in A$ for all subconfigurations $\xi\subset\eta$.

Let $A$ be a hereditary event such that ${\mathbb P}(\mathrm{H}\in A)>0$.
Let $\mathrm{H}_A \sim\mathscr{L}(\mathrm{H}\mvert\mathrm{H}\in
A)$. For instance,
if $A=\{\eta\in\mathfrak{N}\colon|\eta|\le M\}$ for some $M \in
\mathbb{N}$; then $\mathrm{H}_A$ has the same distribution as
$\mathrm{H}$
conditioned on not having more than $M$ points. In many cases $\mathrm{H}_A$
then satisfies $\SSS$, even if the original process $\mathrm{H}$ does not.

The following two lemmas are needed for reducing the problem of
approximating by the process $\mathrm{H}$ to a problem of
approximating by
$\mathrm{H}_A$.
%
%
\begin{lemma}
\label{lemmacondH}
The process $\mathrm{H}_A$ has hereditary density $u_A(\xi)=u(\xi
)\mathbh{1}\{
\xi\in A\}/{\mathbb P}(\mathrm{H}\in A)$ with respect to $\Po_1$ and
conditional intensity $\lambda_A(x\mvert\xi)=\break \lambda(x\mvert\xi
)\mathbh{1}\{\xi+\delta_x\in A\}$, where $u$ and $\lambda$ denote the
density and conditional intensity of $\mathrm{H}$, respectively.
\end{lemma}
\begin{pf}
Note that for all measurable $f:\mathfrak{N}\to\mathbb{R}_+$,
\begin{eqnarray*}
{\mathbb E}f(\mathrm{H}_A) &=&{\mathbb E}\bigl( f(\mathrm{H})
\mvert A\bigr)
\\
&=&\frac{{\mathbb E}( f(\mathrm{H}) \mathbh{1}\{\mathrm{H}\in A\}
)}{{\mathbb P}(\mathrm{H}\in A)}
\\
&=&\int_\mathfrak{N}f(\xi)\frac{u(\xi)\mathbh{1}\{\xi\in A\}
}{{\mathbb
P}(\mathrm{H}\in A)} \Po_1(d
\xi).
\end{eqnarray*}
Furthermore, by the definition of the conditional intensity,
equation~(\ref{eqcif}),
\begin{eqnarray*}
\lambda_A(x\mvert\xi)&=&\frac{u_A(\xi+\delta_x)}{u_A(\xi)}=\frac
{u(\xi+\delta_x) \mathbh{1}\{\xi+\delta_x\in A\}}{u(\xi)
\mathbh{1}\{\xi
\in A\}} =
\lambda(x\mvert\xi)\mathbh{1}\{\xi+\delta_x\in A\},
\end{eqnarray*}
where the last equality follows by the hereditarity of $A$.
\end{pf}

%
\begin{proposition}
\label{propHA}
Let $A$ be a hereditary event, and let $\mathrm{H}_A \sim\mathscr
{L}(\mathrm{H}
\mvert A)$. Then
%
%
\begin{equation}
\label{eqHA} d_{\mathrm{TV}}\bigl(\mathscr{L}(\Xi),\mathscr{L}(\mathrm{H})
\bigr) \le d_{\mathrm
{TV}}\bigl(\mathscr{L}(\Xi),\mathscr{L}(
\mathrm{H}_A)\bigr)+{\mathbb P}(\mathrm{H}\notin A).
\end{equation}
\end{proposition}
\begin{pf}
Note that
\begin{eqnarray*}
&& d_{\mathrm{TV}}\bigl(\mathscr{L}(\mathrm{H}_A),\mathscr{L}(
\mathrm{H})\bigr)
\\
&&\qquad =\sup_{B \in\mathcal{N}}\bigl|{\mathbb P}(\mathrm{H}\in B \mvert
\mathrm{H}\in
A)-{\mathbb P}(\mathrm{H}\in B)\bigr|
\\
&&\qquad =\sup_{B \in\mathcal{N}}\biggl\vert\frac{{\mathbb P}(\mathrm
{H}\in B,
\mathrm{H}\in A)}{{\mathbb P}(\mathrm{H}\in A)}-{\mathbb P}(
\mathrm{H}\in B, \mathrm{H} \in A)- {\mathbb P}(\mathrm{H}\in B,
\mathrm{H}
\notin A)\biggr\vert
\\
&&\qquad =\sup_{B \in\mathcal{N}}\biggl\vert\bigl(1-{\mathbb P}(\mathrm
{H}\in
A)\bigr)\frac
{{\mathbb P}(\mathrm{H}\in B, \mathrm{H}\in A)}{{\mathbb
P}(\mathrm{H}\in
A)}-{\mathbb P}(\mathrm{H}\in B, \mathrm{H}\notin A)
\biggr\vert
\\
&&\qquad \le \max\biggl( \sup_{B \in\mathcal{N}}{\mathbb P}(\mathrm
{H}\notin A)
\frac{{\mathbb P}(\mathrm{H}\in B, \mathrm{H}\in A)}{{\mathbb
P}(\mathrm{H}\in
A)}, \sup_{B \in\mathcal{N}}{\mathbb P}(\mathrm{H}\in B,
\mathrm{H}\notin A) \biggr)
\\
&&\qquad =  {\mathbb P}(\mathrm{H}\notin A)
\end{eqnarray*}
and the triangle inequality yields the claim.
\end{pf}

%
\begin{corollary}
\label{corHA}
For hereditary events $A$ and $A'$ we get
\[
d_{\mathrm{TV}}\bigl(\mathscr{L}(\Xi),\mathscr{L}(\mathrm
{H})\bigr) \le
d_{\mathrm{TV}}\bigl(\mathscr{L}(\Xi_{A}),\mathscr{L}(
\mathrm{H}_{A'})\bigr)+{\mathbb P}(\Xi\notin A)+{\mathbb P}\bigl(
\mathrm{H}\notin A'\bigr).
\]
\end{corollary}

\subsection{General pairwise interaction processes (PIP)}
\label{ssecg-pip}

For PIPs the following hereditary event is very useful. Let $k \in
\mathbb{N}$, $\delta> 0$ and
%
%
\begin{equation}
\label{eqAk} A_k=\Bigl\{\xi\in\mathfrak{N}\colon\sup
_{y\in\mathcal{X}}\xi\bigl(\mathbb{B} (y,\delta/2)\bigr)\le
k\Bigr\},
\end{equation}
that is, we require that
the PIP has at most $k$ points inside any closed ball with radius
$\delta/2$. If $k=1$, this is equivalent to the event that the PIP has
a hard core radius $\delta$.

%
\begin{lemma}
\label{lemmaPAk}
Suppose that $\mathrm{H}\sim\pip(\beta,\varphi)$ satisfies the
conditions $\RS$, $\UB$, $\RC$ and $\IR$ with the constants $C$, $\delta$,
$\gamma$, $r$ and $R$. Then
%
%
\begin{eqnarray}
\label{eqPAk} {\mathbb P}(\mathrm{H}\notin A_k) &=& {\mathbb P} \bigl(
\exists y\in\mathcal{X}\colon\mathrm{H}\bigl(\mathbb{B}(y,\delta /2)\bigr)\ge k+1 \bigr)
\nonumber\\[-8pt]\\[-8pt]
&\le& \frac{\gamma^{(k(k+1))/2}B_\delta^k}{(k+1)!
C^k}{\mathbb E} \bigl(|\mathrm{H}|C^{k|\mathrm{H}|}\bigr),\nonumber
\end{eqnarray}
where $B_\delta=\sup_{y\in\mathcal{X}}\int_{\mathbb{B}(y,\delta
)}\beta
(x) \bolds{\alpha}(dx)$.
\end{lemma}
\begin{pf}
Note that by equation~(\ref{eqpkf})
%
%
\begin{eqnarray}
\label{equglyIntPAk}
&& {\mathbb P}\bigl(\exists y\in\mathcal
{X}\colon\mathrm{H}\bigl(
\mathbb{B}(y,\delta/2)\bigr)\ge k+1\bigr)\nonumber
\\
&&\qquad = e^{-\bolds{\alpha}(\mathcal{X})} \sum_{n=k+1}^\infty
\frac{1}{n!} \int_\mathcal{X}\cdots\int
_\mathcal{X}g_{\delta,k+1}(x_1,
\ldots,x_n)
\\
&&\hspace*{145pt}{}\times u\bigl(\{x_1,\ldots,x_n\}\bigr) \bolds{\alpha}(dx_1)\cdots\bolds{\alpha} (dx_n),
\nonumber
\end{eqnarray}
where
%
%
\begin{eqnarray}\label{eqind}
&& g_{\delta,k+1}(x_1,\ldots,x_n)\nonumber
\\
&&\qquad = \mathbh{1}\bigl\{ \exists y\in\mathcal{X}, \exists\{i_1,
\ldots,i_{k+1}\}\subset\{ 1,\ldots,n\} \colon x_{i_1},
\ldots,x_{i_{k+1}}\in\mathbb{B}(y,\delta/2)\bigr\}
\\
&&\qquad \leq\sum_{\{i_1,\ldots,i_{k+1}\}\subset\{1,\ldots,n\}}\mathbh{1}\bigl\{
\exists y\in\mathcal{X}, x_{i_1},\ldots,x_{i_{k+1}}\in
\mathbb{B}(y,\delta/2)\bigr\}.\nonumber
\end{eqnarray}
For any permutation $(i_1,\ldots,i_n)$ of $(1,\ldots,n)$ the density
can be rewritten as
\begin{eqnarray*}
u\bigl(\{x_1,\ldots,x_n\}\bigr) &=& \Biggl( \prod
_{j=1}^k\beta(x_{i_j})  \Biggr) \biggl(
\prod_{1\le j<l\le k+1} \varphi(x_{i_j},x_{i_l})
\biggr)
\\
&&{} \times u\bigl(\{x_{i_{k+1}},\ldots,x_{i_n}\}\bigr) \prod
_{j=1}^k\prod_{l=k+2}^n
\varphi(x_{i_j},x_{i_l}).
\end{eqnarray*}
Thus by $\RC$ and $\UB$ we obtain
\begin{eqnarray*}
&& u\bigl(\{x_1,\ldots,x_n\}\bigr) \mathbh{1}\bigl\{
\exists y\in\mathcal{X}, x_{i_1},\ldots,x_{i_{k+1}}\in
\mathbb{B}(y,\delta/2)\bigr\}
\\
&&\qquad \le\Biggl(\prod_{j=1}^k
\beta(x_{i_j}) \Biggr) \gamma^{{k+1\choose2}}u\bigl(\{x_{i_{k+1}},
\ldots,x_{i_n}\}\bigr)C^{k(n-k-1)}
\\
&&\quad\qquad{}\times \mathbh{1}\bigl\{\exists y\in\mathcal{X}, x_{i_1},
\ldots,x_{i_{k+1}}\in\mathbb{B}(y,\delta/2)\bigr\}
\\
&&\qquad \le\Biggl(\prod_{j=1}^k
\beta(x_{i_j}) \Biggr) \gamma^{{k+1\choose
2}}C^{k(n-k-1)} u\bigl(
\{x_{i_{k+1}},\ldots,x_{i_n}\}\bigr)
\\
&&\quad\qquad{} \times\mathbh{1}\bigl\{x_{i_1},\ldots,x_{i_k}\in
\mathbb{B} (x_{i_{k+1}},\delta)\bigr\},
\end{eqnarray*}
where the last line follows by the triangle inequality.
Thus in total from equation~(\ref{equglyIntPAk})
\begin{eqnarray*}
\hspace*{-3pt}&& {\mathbb P}(\mathrm{H}\notin A_k)
\\
\hspace*{-3pt}&&\quad \le e^{-\bolds{\alpha}(\mathcal{X})}\sum_{n=k+1}^\infty
\frac{1}{n!}
\\
\hspace*{-3pt}&&\qquad{}\times \sum_{\{i_1,\ldots,i_{k+1}\}\subset\{1,\ldots,n\} }
\underbrace{
\int_\mathcal{X}\cdots\int_\mathcal{X}}_{n-k}
\underbrace{\int_{\mathbb{B}
(x_{i_{k+1}},\delta)}\cdots\int_{\mathbb{B}(x_{i_{k+1}},\delta)}}_{k}
\Biggl(\prod_{j=1}^k\beta(x_{i_j})
\Biggr) \gamma^{{k+1\choose
2}}
\\
\hspace*{-3pt}&&\qquad{} \times C^{k(n-k-1)} u\bigl(\{x_{i_{k+1}},\ldots,x_{i_n}
\}\bigr) \bolds{\alpha}(dx_{i_1})\cdots\bolds{\alpha}(dx_{i_k})
\bolds{\alpha} (dx_{i_{k+1}})\cdots\bolds{\alpha}(dx_{i_n})
\\
\hspace*{-3pt}&&\quad \le\frac{\gamma^{(k(k+1))/2}B_\delta^k}{(k+1)! C^k}
e^{-\bolds{\alpha}(\mathcal{X})}
\\
\hspace*{-3pt}&&\qquad{}\times  \sum_{n=k+1}^\infty
\frac
{n-k}{(n-k)!}C^{k(n-k)}
\underbrace{\int_\mathcal{X}\cdots\int
_\mathcal{X}}_{n-k}u\bigl(\{ x_{k+1},
\ldots,x_n\}\bigr) \bolds{\alpha}(dx_{k+1})\cdots\bolds{
\alpha}(dx_n)
\\
\hspace*{-3pt}&&\quad =\frac{\gamma^{(k(k+1))/2}B_\delta^k}{(k+1)! C^k} {\mathbb
E}\bigl(|\mathrm{H} |C^{k|\mathrm{H}|}\bigr),
\end{eqnarray*}
by equation~(\ref{eqpkf}).
\end{pf}

%
\begin{lemma}
\label{lemmaMk}
Consider the process $\mathrm{H}_{A_k}$, where $\mathrm{H}\sim\pip
(\beta,\varphi)$ satisfies the conditions $\RS$, $\UB$, $\RC$ and $\IR$ with the
constants $C$, $\delta$, $\gamma$, $r$ and $R$. Define
$M_k=C^{m_k}<\infty$ with
\[
m_k=\sup_{x\in\mathcal{X}, \xi\in A_k}\xi\bigl(A(x,r,R)\bigr),
\]
where $A(x,r,R)=\mathbb{B}(x,R)\setminus\mathbb{B}(x,r)$. Then
\[
\lambda_{A_k}(x \mvert\xi)=\lambda(x\mvert\xi)\mathbh{1}\{\xi
+\delta
_x\in A_k\} \le\beta(x)M_k.
\]

This means that the new process $\mathrm{H}_{A_k}$ is locally stable and
hence satisfies condition $\SSS$.

In the Euclidean setting $m_k\le mk$, where
\[
m=\alpha_DD^{D/2} \biggl( \biggl(\frac{R}{\delta}+1
\biggr)^D- \biggl(\frac{r}{\delta}-1 \biggr)^D \biggr).
\]
\end{lemma}

\begin{pf}
By Lemma~\ref{lemmacondH} and conditions $\UB$ and $\IR$ we see that
$\lambda_{A_k}(x\mvert\xi)$ can be bounded by $\beta(x)C^{\xi
(A(x,r,R))}\mathbh{1}\{\xi+\delta_x\in A_k\}\le\beta(x) C^{m_k}$, and
$m_k$ is finite, since it can be bounded by $k$ times the minimal
number of balls with radius $\delta/2$ needed to cover $A(x,r,R)$.

In the Euclidean case, consider a partition $\{Q_i\}_{i=1}^N$ of
$\mathcal{X}$ by cubes of edge length $\delta/\sqrt{D}$. Since the
diameter of each cube is $\delta$, one can cover each cube by a ball
with radius $\delta/2$. Furthermore, a cube can intersect $A(x,r,R)$ if
and only if it is contained in $A(x,r-\delta,R+\delta)$. Thus the
number of cubes intersecting $A(x,r,R)$ can be bounded by the volume of
$A(x,r-\delta,R+\delta)$ divided by the volume of a cube, that is,
\[
\sup_{x\in\mathcal{X}}\frac{|A(x,r-\delta,R+\delta)|}{\delta
^D/D^{D/2}}\le\alpha_DD^{D/2}
\biggl( \biggl(\frac{R}{\delta}+1 \biggr)^D- \biggl(
\frac{r}{\delta}-1 \biggr)^D \biggr)=m.
\]
Thus $m_k\le mk$.
\end{pf}

The next theorem is a generalization of Theorem~\ref{cordtv-pip} that
includes noninhibitory PIPs.

%
\begin{theorem}
\label{thmGeneralPIP}
Assume that $\Xi\sim\pip(\beta,\varphi_1)$ and $\mathrm{H}\sim
\pip
(\beta,\varphi_2)$. Furthermore, assume that they satisfy the
conditions $\RS$, $\UB$, $\RC$ and $\IR$ with the same constants $C$, $\delta
$, $\gamma$, $r$ and $R$. Let $\nu_{A_k}(y) = \mathbb{E}(\nu_{A_k}(y
\mvert\Xi_{A_k}))$ denote the intensity of $\Xi_{A_k}$. Then we have
for any $k \in\mathbb{N}$,
%
%
\begin{eqnarray}
\label{eqGeneralPIP} && d_{\mathrm{TV}}\bigl(\mathscr{L}(\Xi
),\mathscr{L}(\mathrm{H})
\bigr)\nonumber
\\
&&\qquad \le c_1(\lambda)M_k \int_{\mathcal{X}} \int
_\mathcal{X}\beta(x) \nu_{A_k}(y) \bigl|\varphi
_1(x,y)-\varphi_2(x,y)\bigr| \bolds{\alpha}(dx) \bolds{
\alpha}(dy)
\\
&&\quad\qquad{}+\frac{\gamma^{(k(k+1))/2}B_\delta
^k}{(k+1)! C^k}{\mathbb E} \bigl(|\Xi|C^{k|\Xi|}+|
\mathrm{H}|C^{k|\mathrm{H}|} \bigr),\nonumber
\end{eqnarray}
where $M_k$ is defined in Lemma~\ref{lemmaMk}, and $c_1(\lambda)$ is
given in inequality~(\ref{eqc1}) with
\[
c \leq M_k\int_\mathcal{X}\beta(x) \bolds{
\alpha}(dx)
\]
and
\[
\varepsilon\leq M_k\sup_{y\in\mathcal{X}} \biggl( \int
_{\mathcal
{X}\setminus\mathbb{B}(y,\delta)} \beta(x)\bigl|\varphi_2(x,y)-1\bigr|
\bolds{
\alpha}(dx)+\int_{\mathbb
{B}(y,\delta)}\beta(x) \bolds{\alpha}(dx) \biggr).
\]
\end{theorem}

Note that $k \in\mathbb{N}$ can be chosen such that the two terms in
(\ref{eqGeneralPIP}) are best balanced. There is also some freedom in
the choice of $\delta$ in $\IR$. In particular it is always possible to
choose a lower $\delta$ at no cost for $\gamma$, that is, the last
term in (\ref{eqGeneralPIP}) can be made arbitrarily small by letting
$\delta\to0$, which on the other hand leads to an explosion of $M_k$.

%
\begin{remark}
Usually the expectations in (\ref{eqGeneralPIP}) are not easy to
compute, but at least Ruelle stability guarantees their finiteness. Let
$u$ be the density of $\Xi$, let $\psi^*$ be as in $\RS$, and write
$\bolds{\alpha}(\psi^*)=\int_\mathcal{X}\psi^*(x) \bolds{\alpha
}(dx)$. Then
for a constant $c^{**}>0$
\begin{eqnarray*}
{\mathbb E}\bigl(|\Xi|C^{k|\Xi|}\bigr)&\le& c^{**}e^{-\bolds{\alpha
}(\mathcal
{X})}
\sum_{n=0}^\infty\frac{nC^{kn}\bolds{\alpha}(\psi^*)^n}{n!}
\\
&=& c^{**}C^k\bolds{\alpha}\bigl(\psi^*\bigr)e^{C^k\bolds{\alpha
}(\psi
^*)-\bolds{\alpha}(\mathcal
{X})}
<\infty.
\end{eqnarray*}
\end{remark}

%
\begin{remark}
If one is interested in a very specific Gibbs processes model, then the
estimates in Lemmas~\ref{lemmaPAk} and~\ref{lemmaMk} may be
improved. See, for instance, Section~\ref{sseclenny-jones}, where we
treat Lennard--Jones type processes.
\end{remark}

By a slight adaptation in the proof of Theorem~\ref{cordtv-pip}, we
can get a nicer result for hard core PIPs.
%
%
\begin{theorem}
\label{thmPIPHC}
Let $\Xi\sim\pip(\beta_1,\varphi_1)$ and $\mathrm{H}\sim\pip
(\beta
_2,\varphi_2)$. Assume that both processes have a
hard core radius of $\delta>0$, and satisfy $\UB$ and $\IR$. Then
%
%
\begin{eqnarray}
\label{eqpiphc}
&& d_{\mathrm{TV}}\bigl(\mathscr{L}(\Xi),\mathscr{L}(\mathrm{H})\bigr)
\nonumber\\[-8pt]\\[-8pt]
&&\qquad \le c_1(\lambda) M_1 \int
_\mathcal{X}\int_{\mathcal{X}} \beta(x) \nu(y) \bigl|
\varphi_1(x,y)-\varphi_2(x,y)\bigr| \bolds{\alpha}(dx) \bolds{
\alpha}(dy).\nonumber
\end{eqnarray}
\end{theorem}

\begin{pf*}{Proof of Theorem~\ref{thmGeneralPIP}}
We adapt the proof of Theorem~\ref{cordtv-pip} to compute a bound for
$d_{\mathrm{TV}}(\mathscr{L}(\Xi_{A_k}),\mathscr{L}(\mathrm
{H}_{A_k}))$. Let $\xi
=\sum_{i=1}^n\delta_{y_i}$. Since both $\varphi_1$ and $\varphi_2$
are bounded by the same~$C$,
\begin{eqnarray*}
&& \Biggl(\prod_{i=1}^{j-1}
\varphi_1(x,y_i) \Biggr) \Biggl( \prod
_{i=j+1}^n\varphi_2(x,y_i)
\Biggr) \mathbh{1}\{\xi+\delta_x \in A_k\}
\\
&&\qquad \le C^{\xi(\mathbb{B}(x,R)\setminus\mathbb{B}(x,r))}\mathbh{1}\{
\xi+\delta_x \in A_k
\} \le C^{m_k}= M_k.
\end{eqnarray*}
Thus
%
%
\begin{eqnarray}
\label{eqcondpips}
&& d_{\mathrm{TV}}\bigl(\mathscr{L}(\Xi_{A_k}),\mathscr{L}(\mathrm{H}_{A_k})\bigr)
\nonumber\\[-8pt]\\[-8pt]
&&\qquad \leq c_1(\lambda)M_k \int
_{\mathcal{X}} \int_\mathcal{X}\beta(x)
\nu_{A_k}(y) \bigl|\varphi_1(x,y)-\varphi_2(x,y)\bigr|
\bolds{\alpha}(dx) \bolds{\alpha}(dy).\nonumber
\end{eqnarray}

Since by Lemma~\ref{lemmaMk} $|\lambda_{A_k}(x\mid\xi)-\lambda
_{A_k}(x\mid\eta)| \le\beta(x)M_k$ for all $x\in\mathcal{X}$ and
\mbox{$\xi,\eta\in\mathfrak{N}$,} we have $c \leq M_k\int_\mathcal
{X}\beta(x) \bolds{\alpha}(dx)$. By Remark~\ref{remarkesssup} the
supremum in the formula for $\varepsilon$ can be replaced by an
essential supremum
with respect to $\mathscr{L}(\Xi_{A_k}) + \mathscr{L}(\mathrm{H}_{A_k})$.
This implies that it is enough to take the supremum only over $\xi,\eta\in A_k$. Note that if $d(x,y)>\delta$, $\xi+\delta_x\in A_k$
and $\xi+\delta_y\in A_k$, then also $\xi+\delta_x+\delta_y\in
A_k$. Therefore, by using $\lambda(x\mid\xi+\delta_y)=\lambda
(x\mid\xi)\varphi_2(x,y)$, we obtain
\begin{eqnarray*}
\varepsilon&=& \sup_{\xi,\eta\in A_k, \llVert \xi-\eta\rrVert
=1}\int_\mathcal{X}\bigl|
\lambda_{A_k}(x\mid\eta)-\lambda_{A_k}(x\mid\xi)\bigr| \bolds{
\alpha}(dx)
\\
&=& \sup_{y\in\mathcal{X}, \xi+\delta_y\in A_k} \int_\mathcal{X}
\bigl| \lambda(x
\mid\xi+\delta_y)\mathbh{1}\{\xi+\delta_x+\delta
_y\in A_k\}
\\
&&\hspace*{104pt}{} -\lambda(x\mid\xi)\mathbh{1}\{ \xi+\delta_x\in A_k
\} \bigr| \bolds{\alpha}(dx)
\\
&\le& \sup_{y\in\mathcal{X}, \xi+\delta_y\in A_k} \biggl(\int
_{\mathcal{X}\setminus\mathbb{B}(y,\delta)}
\lambda(x\mid\xi)\bigl|\varphi_2(x,y)-1\bigr|\mathbh{1}\{\xi+
\delta_x+\delta_y\in A_k\} \bolds{\alpha
}(dx)
\\
&&\hspace*{61pt}+\int_{\mathbb{B}(y,\delta)} \bigl| \lambda(x\mid\xi+\delta_y)
\mathbh{1}\{\xi+\delta_x+\delta_y\in A_k\}
\\
&&\hspace*{171pt}{} -\lambda(x\mid\xi)\mathbh{1}\{\xi+\delta_x\in
A_k\} \bigr| \bolds{\alpha}(dx) \biggr)
\\
&\le& M_k\sup_{y\in\mathcal{X}} \biggl( \int
_{\mathcal
{X}\setminus\mathbb{B}(y,\delta)} \beta(x)\bigl|\varphi_2(x,y)-1\bigr|
\bolds{
\alpha}(dx)+\int_{\mathbb
{B}(y,\delta)}\beta(x) \bolds{\alpha}(dx) \biggr).
\end{eqnarray*}

The claim now follows by applying Corollary~\ref{corHA} with
$A=A'=A_k$ and Lemma~\ref{lemmaPAk}.
\end{pf*}

\begin{pf*}{Proof of Theorem~\ref{thmPIPHC}}
Since $\mathscr{L}(\Xi)=\mathscr{L}(\Xi_{A_1})$ and $\mathscr
{L}(\mathrm{H})=\mathscr{L}(\mathrm{H}_{A_1})$, the statement
follows from
inequality~(\ref{eqcondpips}).
\end{pf*}

%
\begin{example}
\label{exmsstrauss}
Let $\mathcal{X}\subset\mathbb{R}^D$. A PIP is called a \emph
{multi-scale Strauss process}, if its interaction function $\varphi
(x,y)$ depends only on $\llVert x-y \rrVert$, is piecewise constant
and takes
only finitely many values. We restrict ourselves to bi-scale Strauss
processes, that is, the interaction function is given by
\[
\varphi(x,y)= %
\cases{\gamma, &\quad if $\llVert x-y \rrVert\le r$,
\cr
C, &\quad if $r < \llVert x-y \rrVert\le R$,
\cr
1, &\quad if $\llVert x-y
\rrVert> R$}
\]
for some constants $0\le\gamma\le1$, $C\ge0$ and $0<r<R$. To ensure
$\RS$, we furthermore require that $C\le\gamma^{-1/(2m)}$ with
\[
m=m(r,R;D)=\alpha_DD^{D/2} \biggl(\frac{R}{r}+1
\biggr)^D.
\]

Note that $m$ is the same as in Lemma~\ref{lemmaMk} with $\delta=r$.
$\RS$ then follows by a criterion of \citet{kpr2012}, Section~2.3. The
authors use $m$ as a bound on the maximal number of cubes with edge
length $\delta/\sqrt{D}$ that intersect the annulus $A(0,r-\delta,R+\delta)$, as we did in Lemma~\ref{lemmaMk}.

To illustrate Theorem~\ref{thmGeneralPIP}, let $\Xi$ and $\mathrm
{H}$ be
bi-scale Strauss processes with constant $\beta$, the same $\gamma$,
$r$, $R$, and with $1 \leq C_\mathrm{H}\leq C_\Xi$. The ingredients for
computing $c_1(\lambda)$ are
\[
\varepsilon=\alpha_D\beta C_\mathrm{H}^{mk}
\bigl(r^D+(C_\mathrm{H} -1) \bigl(R^D-r^D
\bigr) \bigr) \quad\mbox{and}\quad c=\beta C_\mathrm{H} ^{mk}|
\mathcal{X}|.
\]
Since $\gamma\le C_\Xi^{-2m}$, Theorem~\ref{thmGeneralPIP} yields
\begin{eqnarray*}
d_{\mathrm{TV}}\bigl(\mathscr{L}(\Xi),\mathscr{L}(\mathrm
{H})\bigr) &\le&
c_1(\lambda)\alpha_D\beta C_{\Xi}^{mk}{
\mathbb E}\bigl(|\Xi_{A_k}|\bigr) (C_\Xi-C_\mathrm{H})
\bigl(R^D-r^D\bigr)
\\
&&{}+\frac{\alpha_D\beta^kr^{Dk}}{(k+1)!}C_\Xi
^{-mk^2-(m+1)k}{\mathbb E} \bigl(|
\Xi|C_\Xi^{k|\Xi|}+|\mathrm{H}|C_\Xi^{k|\mathrm
{H}|}
\bigr).
\end{eqnarray*}
\end{example}

\section{Applications}
\label{secapplications}

\subsection{Lennard--Jones type processes}
\label{sseclenny-jones}

In this subsection let $(\mathcal{X},\bolds{\alpha})$ be a compact
subset of
$\mathbb{R}^D$ equipped with Lebesgue measure. We say a PIP is of
\emph{Lennard--Jones type} [see \citet{ruelle69}], if its interaction
function can be written as $\varphi(x,y)=\exp(-bV(\llVert x-y
\rrVert))$, and
the pair potential $V$ satisfies the following conditions:
\begin{longlist}[(2)]
\item[(1)] There exist $r\leq R$ and a $\varrho> D$ such that
\begin{eqnarray*}
V(x) &\ge&\llVert x \rrVert^{-\varrho} \qquad\mbox{for } \llVert x \rrVert
\le r,
\\
V(x) &\ge& -\llVert x \rrVert^{-\varrho} \qquad\mbox{for } \llVert x \rrVert
\ge R.
\end{eqnarray*}
\item[(2)] $V(x)\ge-M$ for a $M\ge0$ and for all $x\ge0$, that is,
the interaction function~$\varphi$ is bounded from above by $e^{bM}$.
\end{longlist}

The technique used in Section~\ref{ssecg-pip} stands and falls
with a good estimate on the term $\sup_{\xi\in\mathfrak{N}}\lambda
_{A_k}(x \mvert\xi)$. The next lemma gives a neat replacement of
Lemma~\ref{lemmaMk} for Lennard--Jones type processes. 

%
\begin{lemma}
\label{lemmaLJMk}
Assume that $\mathrm{H}$ is a $\pip$ of Lennard--Jones type with constants
$\varrho,r,R,M$. Choose a positive $\delta\leq r$ such that also
$\delta<R/2$. Define
%
%
\begin{equation}
\label{eqLJMk} M_k = \exp\biggl(bk \biggl(mM+\frac{\alpha
_DD}{\varrho-D}
\biggl(\frac{\sqrt{D}}{\delta} \biggr)^D\frac{(R-\delta
)^{D-1}}{(R-2\delta)^{\varrho-1}} \biggr) \biggr),
\end{equation}
where
\[
m=m(r,R;\delta)=\alpha_DD^{D/2} \biggl( \biggl(
\frac{R}{\delta
}+1 \biggr)^D- \biggl(\frac{r}{\delta}-1
\biggr)^D \biggr)\mathbh{1}\{ r<R\}.
\]
Then for all $x\in\mathcal{X}$ and for all $\xi\in\mathfrak{N}$ we have
\[
\lambda_{A_k}(x\mid\xi)\le\beta(x)M_k.
\]
\end{lemma}

\begin{pf}
Since $\varphi\le e^{bM}$, we obtain analogously as in the proof of
Lemma~\ref{lemmaMk} and by using the translation invariance of $V$ that
%
%
\begin{eqnarray}\label{eqprooflenny1}
\qquad \sup_{\xi\in A_k}\lambda(x\mid\xi)&\le&\beta(x)\sup
_{\xi\in A_k} \exp\biggl(bM\xi\bigl(A(x,r,R)\bigr)-b \sum
_{y\in\xi, \llVert y
\rrVert\ge R } V\bigl(\llVert y \rrVert\bigr) \biggr)
\nonumber\\[-8pt]\\[-8pt]
&\le&\beta(x)e^{bMmk}\sup_{\xi\in A_k} \exp\biggl(-b \sum
_{y\in
\xi, \llVert y \rrVert\ge R } V\bigl(\llVert y \rrVert\bigr)
\biggr).\nonumber
\end{eqnarray}

Let $\{Q_z\}_{z\in\mathbb{Z}^D}$ denote the partition of $\mathbb
{R}^D$ into cubes of edge length $\delta/\sqrt{D}$ and centre points
$\delta/\sqrt{D} \mathbb{Z}^D$. Since $\xi\in A_k$, each cube
contains at most $k$ points. A cube intersects $\mathbb{B}(0,R)^c$ if and
only if its center point is contained in $\mathbb{B}(0,R-\delta
/2)^c$. For
$x=(x_1,\ldots,x_D)\in\mathbb{R}^D$, denote $\llVert x \rrVert
_{\max}=\max_{i=1,\ldots,D}\vert x_i \vert$. Thus
%
%
\begin{eqnarray}\label{eqlennyjonesMk}
&& \sup_{\xi\in A_k}\sum_{y\in\xi, \llVert y \rrVert\ge R }
-V\bigl(\llVert y \rrVert\bigr)\nonumber
\\
&&\qquad \le \mathop{\sum_{z\in
\mathbb{Z}^D}}_{\llVert z \rrVert \ge\sqrt{D}(R/\delta
-1/2)}k\sup_{\llVert x-z \rrVert_{\max} \le 1/2} \biggl(\frac{\delta}{\sqrt{D}}\llVert x \rrVert\biggr)^{-\varrho}
\\
&&\qquad =k \biggl(\frac{\delta}{\sqrt{D}} \biggr)^{-\varrho} \mathop{\sum
_{z\in\mathbb{Z}^D}}_{\llVert z \rrVert
\ge\sqrt {D}(R/\delta-1/2)}\sup_{\llVert x-z \rrVert_{\max} \le1/2}\llVert x \rrVert
^{-\varrho}.\nonumber
\end{eqnarray}

Consider the function $g(x)=(\llVert x \rrVert-\sqrt{D})^{-\varrho
}$. Since
\[
\inf_{\llVert x-z \rrVert_{\max} \le 1/2}g(x) \ge\biggl
(\llVert z \rrVert-
\frac{\sqrt{D}}{2} \biggr)^{-\varrho} \ge\sup_{\llVert x-z
\rrVert_{\max
} \le 1/2} \llVert
x \rrVert^{-\varrho},
\]
the last sum in (\ref{eqlennyjonesMk}) can be bounded by an integral
over the function $g$. For any $a>0$ we get
\begin{eqnarray*}
&& \int_{\mathbb{B}(0,a)^c}\bigl(\llVert x \rrVert-\sqrt{D}
\bigr)^{-\varrho} \,dx
\\
&&\qquad = \alpha_DD\int_a^\infty(r-
\sqrt{D})^{-\varrho}r^{D-1} \,dr
= \alpha_DD\int_{a-\sqrt{D}}^\infty\biggl(
\frac{r+\sqrt{D}}{r} \biggr)^{D-1}r^{D-1-\varrho} \,dr
\\
&&\qquad \le \alpha_DD \biggl(\frac{a}{a-\sqrt{D}} \biggr)^{D-1}\int
_{a-\sqrt{D}}^\infty r^{D-\varrho-1} \,dr
\\
&&\qquad = \alpha_DD \biggl(\frac{a}{a-\sqrt{D}} \biggr)^{D-1}
\frac{1}{\varrho-D}\frac{1}{(a-\sqrt{D})^{\varrho-D}}.
\end{eqnarray*}
In order to catch all cubes in (\ref{eqlennyjonesMk}), the
integration must begin at $a=\sqrt{D}(R/\delta-1)$, which together
with (\ref{eqprooflenny1}) yields the claim.
\end{pf}

We may also give a more explicit bound on $\mathbb{P}(\mathrm
{H}\notin
A_k)$ than inequality~(\ref{eqPAk}) in the case of a Lennard--Jones
type process.

%
\begin{lemma}
\label{lemmaLJ-PAk}
Assume that $\mathrm{H}$ is a $\pip$ of Lennard--Jones type with constants
$\varrho,r,R,M$. Then for any $\delta<\min(r,R/2)$
%
%
\begin{equation}
\label{eqLJ-PAk} {\mathbb P}(\mathrm{H}\notin A_k)\le\biggl(\int
_\mathcal{X}\beta(x) \,dx \biggr)\sum_{j=k+1}^\infty
\frac{1}{j!}B_\delta^{j-1} \exp\bigl(-j^2bL(
\delta) \bigr),
\end{equation}
where $B_\delta$ is given in Lemma~\ref{lemmaPAk} and
\begin{eqnarray*}
L(\delta)&=&\frac{1}{4}\delta^{-\varrho}-\delta^{-D} \biggl(M
\alpha_DD^{D/2} \bigl((R+ \delta)^D-(r-
\delta)^D \bigr) \mathbh{1}\{ r<R\}
\\
&&\hspace*{133pt}{}+ \frac{\alpha_DD^{D/2+1}}{\varrho-D}\frac{(R-\delta
)^{D-1}}{(R-2\delta)^{\varrho-1}} \biggr),
\end{eqnarray*}
which is positive for reasonably small $\delta$, since $\varrho>D$.
\end{lemma}
\begin{pf}
Let $\widetilde{A}_k=A_k\setminus A_{k-1}$ for $k\ge2$. The sets
$\widetilde{A}_k$, $k\ge2$, are pairwise disjoint and $\mathfrak
{N}\setminus
A_l= \bigcup_{k=l+1}^\infty\widetilde{A}_k$ for all $l\ge1$. Then by
(\ref{eqpkf})
%
%
\begin{eqnarray}
\label{equglyLJ} \qquad&& {\mathbb P}(\mathrm{H}\in\widetilde {A}_k)\nonumber
\\
&&\qquad =e^{-\vert \mathcal{X}
\vert}\sum_{n=k}^\infty\frac{1}{n!}\int
_\mathcal{X}\cdots\int_\mathcal{X}\mathbh{1}
\bigl\{\{x_1,\dots,x_n\}\in A_k\bigr\}
\\
&&\hspace*{128pt}{}\times g_{\delta,k}(x_1,\ldots,x_n) u
\bigl(\{x_1,\dots,x_n\}\bigr) \,dx_1\cdots
dx_n,\nonumber
\end{eqnarray}
where $g_{\delta,k}$ has been defined in (\ref{eqind}).
For any permutation $(i_1,\dots,i_n)$ of $(1,\dots,n)$ the density
can be rewritten in a similar way as in the proof of Lemma~\ref
{lemmaPAk} as
\begin{eqnarray*}
&& u\bigl(\{x_1,\dots,x_n\}\bigr)
\\
&&\qquad = \biggl(\prod
_{1\le j<l\le k}\varphi(x_{i_j},x_{i_l}) \biggr)
\Biggl(\prod_{j=1}^k
\lambda\bigl(x_{i_j}\mid\{ x_{i_{k+1}},\dots,x_{i_n}\}
\bigr) \Biggr)u\bigl(\{x_{i_{k+1}},\dots,x_{i_n}\}\bigr).
\end{eqnarray*}
Since $\delta\le r$, for all $x,y$ with $\llVert x-y \rrVert\le
\delta$, we
have $\varphi(x,y)\le\exp(-b\delta^{-\varrho})$.
The term $g_{\delta,k}(x_1,\ldots,x_n)$ can be bounded as in (\ref
{eqind}), and by using Lemma~\ref{lemmaLJMk} and the\vadjust{\goodbreak} triangle
inequality, we get
\begin{eqnarray*}
&& {\mathbb P}(\mathrm{H}\in\widetilde{A}_k)
\\[-3pt]
&&\qquad \le e^{-\vert
\mathcal{X} \vert}\sum_{n=k}^\infty
\frac{1}{n!} \sum_{\{i_1,\ldots,i_{k}\}\subset\{
1,\ldots,n\}}\underbrace{\int_\mathcal{X}\cdots\int_\mathcal{X}}_{n-k}
\int_\mathcal{X}\underbrace{\int_{\mathbb{B}(x_{i_k},\delta
)}\cdots
\int_{\mathbb{B}(x_{i_k},\delta)}}_{k-1} e^{-b\delta^{-\varrho
}{k\choose2}}
\\[-4pt]
&&\qquad\quad{} \times \Biggl(\prod_{j=1}^k
\beta(x_{i_j}) \Biggr)M_k^k u\bigl(
\{x_{i_{k+1}},\dots,x_{i_n}\}\bigr) \,dx_{i_1}\cdots
dx_{i_n}
\\
&&\qquad \le\frac{1}{k!} \biggl(\int_\mathcal{X}\beta(x) \,dx
\biggr) B_\delta^{k-1}e^{-b\delta^{-\varrho}{k\choose2}}M_k^k
e^{-\vert\mathcal{X} \vert}
\\[-2pt]
&&\quad\qquad{} \times \sum_{n=k}^\infty
\frac{1}{(n-k)!}\underbrace{\int_\mathcal{X}\cdots\int
_\mathcal{X}}_{n-k}u\bigl(\{x_{k+1},
\dots,x_n\}\bigr) \,dx_{k+1}\cdots dx_n
\\[-2pt]
&&\qquad = \frac{1}{k!} \biggl(\int_\mathcal{X}\beta(x) \,dx
\biggr)B_\delta^{k-1}e^{-b\delta^{-\varrho}{k\choose2}}M_k^k.
\end{eqnarray*}
For $k\ge2$ we have ${k\choose2}=(k^2/2)(1-1/k)\ge k^2/4$. Thus
Lemma~\ref{lemmaLJMk} yields $\exp(-b\delta^{-\varrho
}{k\choose2} )M_k^k\le\exp(-bk^2L(\delta))$. The statement
then follows by ${\mathbb P}(\mathrm{H}\notin A_l)=\sum
_{k=l+1}^\infty
{\mathbb P}(\mathrm{H}\in\widetilde{A}_k)$.
\end{pf}

Putting the pieces together as we did for the proof of Theorem~\ref
{thmGeneralPIP} we obtain the corresponding upper bound. For the sake
of simplicity we formulate the following result for the special case of
the classical Lennard--Jones process in three dimensions, where the
pair potential is of the form
\[
V_R(x)= \biggl(\frac{R}{\llVert x \rrVert} \biggr)^{12}- \biggl(
\frac
{R}{\llVert x \rrVert} \biggr)^6
\]
for $R>0$ and $x \in\mathbb{R}^3$. In this case we can choose $r=R$
and $\varrho=6$.
%
%
\begin{theorem}
Let $\Xi\sim\pip(\beta,\varphi_1)$ and $\mathrm{H}\sim\pip
(\beta,\varphi_2)$ be classical Lennard--Jones processes with interaction
functions $\varphi_i(x) =\exp(-b_iV_{R_i}(x) )$ for
$i=1,2$. Let $\nu_{A_k}(y) = \mathbb{E}(\nu_{A_k}(y \mvert\Xi_{A_k}))$
denote the intensity of $\Xi_{A_k}$. Then we have for any $k \in
\mathbb{N}$ and $\delta< \min(R_1,R_2)/2$ that
%
%
\begin{eqnarray}
\label{eqex-lj}
&& d_{\mathrm{TV}}\bigl(\mathscr{L}(\Xi),\mathscr{L}(\mathrm{H})\bigr)\nonumber
\\[-3pt]
&&\qquad \le c_1(\lambda)\exp\biggl(b_2k4\pi\sqrt{3}
\frac{(R_2-\delta)^2}{\delta^3(R_2-2\delta)^5} \biggr)
\nonumber\\[-10pt]\\[-10pt]
&&\quad\qquad{} \times \int_\mathcal{X}\int_\mathcal{X}
\beta(x)\nu_{A_k}(y)\bigl\vert\varphi_{1}\bigl(\llVert
x-y \rrVert\bigr)-\varphi_{2}\bigl(\llVert x-y \rrVert\bigr)
\bigr\vert \,dx \,dy\nonumber
\\[-3pt]
&&\quad\qquad{}+ \biggl(\int_\mathcal{X}\beta(x) \,dx \biggr)
\sum_{j=k+1}^\infty\frac{1}{j!}B_\delta^{j-1}
\bigl(e^{-j^2b_1L_1(\delta)}+e^{-j^2b_2L_2(\delta)} \bigr),
\nonumber
\end{eqnarray}
where
\[
L_i(\delta)=\frac{1}{4}\delta^{-6}-4\pi\sqrt{3}
\frac{(R_i-\delta
)^2}{\delta^{3}(R_i-2\delta)^5}
\]
for $i=1,2$.
Note that $\nu_{A_k}$ may be bounded in a crude manner by
\[
\nu_{A_k}(y)\le\beta(y)\exp\biggl(b_1k4\pi\sqrt{3}
\frac
{(R_1-\delta)^2}{\delta^3(R_1-2\delta)^5} \biggr).
\]
\end{theorem}

%
\begin{remark}
Typically any endeavors to make the second summand in inequality~(\ref{eqex-lj}) small make the exponential factor in the first summand quite
large so that the bounds are mainly useful in an asymptotic setting
where the interaction functions are very close. Note, however, that if
$R_1$ and $R_2$ are large, so that $\delta$ may be chosen quite a bit
larger than $1$ but still substantially smaller than $R_i$, choosing a
large $k$ results in a situation where the second summand is close to
zero and the exponential factor is close to one.
\end{remark}

\subsection{The hard core process as limit of area interaction processes}
\label{ssecaip}
For simplicity let again $(\mathcal{X},\bolds{\alpha})$ be a compact subset
of $\mathbb{R}^D$ with Lebesgue measure. In this subsection let
$\mathbb{B}
(x,R)$ always denote the closed ball in $\mathbb{R}^D$ rather than in
$\mathcal{X}$. Let $\mathrm{H}$ be a Strauss process with parameters
$R,\beta_0 > 0$ and $\gamma_0=0$ (hard core case). Let furthermore
$\Xi:= \Xi_{\beta,\gamma}$ be an area interaction process with
parameters $R/2,\beta,\gamma$, where $\gamma\in(0,1]$. The
unnormalized density of such a process is given by
\[
\tilde{u}(\xi) = \beta^{\vert \xi\vert} \gamma^{-\vert
\bigcup_{y \in\xi} \mathbb{B}(y,R/2) \vert}
\]
and the conditional intensity is therefore
\[
\nu(x\mvert\xi)=\beta\gamma^{-|\mathbb{B}(x,R/2)\setminus\bigcup
_{y\in
\xi}\mathbb{B}(y,R/2)|};
\]
see \citet{bv95} for more details. The authors show that $\mathscr
{L}(\Xi_{\beta,\gamma}) \to\mathscr{L}(\mathrm{H})$ weakly as
$\beta,
\gamma\to0$ in such a way that $\beta\gamma^{-\alpha_D (R/2)^D}
\to\beta_0$ (it is easily seen that the hard core process referred to
by Baddeley and van Lieshout is in fact the Strauss hard core process
we use). We derive a rate for this convergence.

%
\begin{theorem}
\label{thmAIP-HC}
Let $\Xi$ and $\mathrm{H}$ be as above. Then
\begin{eqnarray*}
\label{eqAIP-HC} && d_{\mathrm{TV}}\bigl(\mathscr{L}(\Xi),\mathscr
{L}(\mathrm{H})
\bigr)
\\
&&\qquad \leq c_1(\lambda) \bigl( \bigl| \beta\gamma^{-\alpha_D(R/2)^D} -
\beta_0 \bigr||\mathcal{X}| + \beta\gamma^{-\alpha_D (R/2)^D} {\mathbb E}|
\Xi| I_D(R,\gamma) \bigr),
\end{eqnarray*}
where
\[
I_D(R,\gamma):= \int_{\mathbb{B}(0,R)}
\gamma^{\vert \mathbb
{B}(x,R/2) \cap\mathbb{B} (0,R/2) \vert} \,dx \leq2 \alpha_DD
R^{D-1}\log\bigl(
\gamma^{-\alpha_D}\bigr)^{-1/D}.
\]
\end{theorem}

\begin{pf}
For the difference between the conditional intensities we obtain
\[
\nu(x \mvert\xi) - \lambda(x \mvert\xi) = %
\cases{ \beta
\gamma^{-\vert \mathbb{B}(x,R/2) \vert} - \beta_0, &\quad if
$\dist(x,\xi) > R$,
\vspace*{3pt}\cr
\beta\gamma^{-\vert \mathbb{B}(x,R/2) \setminus\bigcup_{y \in\xi}
\mathbb{B} (y,R/2) \vert}, &\quad if $\dist(x,\xi) \leq R$,}
\]
where $\dist(x,A) = \inf_{y \in A} \llVert x-y \rrVert$ for any $A
\subset
\mathbb{R}^D$.
Therefore
\begin{eqnarray*}
&& \int\bigl| \nu(x \mvert\xi) - \lambda(x \mvert\xi) \bigr| \,dx
\\
&&\qquad \leq\int\bigl| \beta\gamma^{-\vert \mathbb{B}(x,R/2) \vert} - \beta
_0 \bigr| \,dx
\\
&&\quad\qquad {}+ \beta\gamma^{-\alpha_D (R/2)^D} \int_{\{\tilde{x}:
\dist(\tilde{x},\xi) \leq R\}}
\gamma^{\vert \mathbb{B}(x,R/2) \cap\bigcup_{y \in\xi} \mathbb
{B}(y,R/2) \vert} \,dx.
\end{eqnarray*}

The last integral may be bounded further by
\begin{eqnarray*}
&& \vert\xi\vert\int_{\mathbb{B}(0,R)} \gamma^{\vert \mathbb
{B}(x,R/2) \cap\mathbb{B} (0,R/2) \vert} \,dx
\\
&&\qquad \leq
\vert\xi\vert\int_{\mathbb{B}(0,R)} \gamma^{\alpha_D(R-\llVert
x \rrVert)^D/2^D} \,dx
= \vert\xi\vert\alpha_DD \int_0^R
\gamma^{\alpha_D(R-r)^D/2^D} r^{D-1} \,dr.
\end{eqnarray*}
By the substitution $y=\log(\gamma^{-\alpha_D}) (\frac
{R-r}{2} )^D$ this is equal to
\begin{eqnarray*}
&&\vert\xi\vert\alpha_D D \log\bigl(\gamma^{-\alpha
_D}
\bigr)^{-1/D}\frac
{2}{D}\int_0^{\log(\gamma^{-\alpha_D}) (R/2)^D}
e^{-y}y^{(1/D)-1}
\\
&&\hspace*{192pt}{}\times \biggl(R-2 \biggl(\frac{y}{\log(\gamma^{-\alpha
_D})} \biggr)^{1/D}
\biggr)^{D-1} \,dy
\\
&&\qquad \le2\vert\xi\vert\alpha_DD R^{D-1}
\log\bigl(\gamma^{-\alpha_D}\bigr)^{-1/D}\frac{1}{D}\int_0^\infty e^{-y}y^{(1/D)-1} \,dy
\\
&&\qquad \le2\vert\xi\vert\alpha_DD R^{D-1}
\log\bigl(\gamma^{-\alpha_D}\bigr)^{-1/D}.
\end{eqnarray*}
The last inequality holds since the integral is equal to $\Gamma(1/D)$
and the functional equality of the Gamma function yields $(1/D) \Gamma
(1/D)=\Gamma(1+1/D)\le1$ for all
$D\ge1$.
The result follows now from Theorem~\ref{thmmain2}.
\end{pf}

The next proposition shows that in the case $\beta\gamma^{-\alpha_D
(R/2)^D} = \beta_0$ the above rate is sharp. Define $\mathcal
{X}^{(-R)}:= \{x \in\mathcal{X}\colon\dist(x,\mathcal{X}^c) \geq
R \}$ and choose $R_0$ such that $\vert \mathcal{X}^{(-R_0)} \vert > 0$.
%
%
\begin{proposition}
\label{propAIP-lbound}
Let $\Xi$ and $\mathrm{H}$ be as above. Assume that $\beta\gamma
^{-\alpha
_D (R/2)^D} = \beta_0$ and $R \leq R_0$. Then there exists a positive
constant $\kappa$ such that
%
%
\begin{equation}
\label{eqAIP-lbound} d_{\mathrm{TV}}\bigl(\mathscr{L}(\Xi
),\mathscr{L}(\mathrm{H})
\bigr) \ge\kappa I_D(R,\gamma).
\end{equation}
\end{proposition}

\begin{pf}
Define $A=\{\xi\in\mathfrak{N}\colon\exists\{x,y\}\subset\xi,
\llVert x-y \rrVert\leq R\}$.  Note that\break \mbox{${\mathbb P}(\mathrm{H}\in
A)=0$}. Hence
\begin{eqnarray*}
&& d_{\mathrm{TV}}\bigl(\mathscr{L}(\Xi),\mathscr{L}(\mathrm
{H})\bigr)
\\
&&\qquad =  \sup _{B\in\mathcal
{N}}\bigl|{\mathbb P}(\Xi\in B)- {\mathbb P}(\mathrm{H}\in B)\bigr|
\ge \bigl|{\mathbb P}(\Xi\in A)- {\mathbb P}(\mathrm{H}\in A)\bigr|={\mathbb
P}(\Xi\in A).
\end{eqnarray*}
Denote by $c_\Xi$ the normalizing constant of the density of $\Xi$. Then
\begin{eqnarray*}
&& {\mathbb P}(\Xi \in A)
\\
&&\qquad \ge{\mathbb P}\bigl(\Xi\in A, \vert\Xi\vert=2\bigr)
\\
&&\qquad = c_{\Xi} \frac{e^{-\vert \mathcal{X} \vert}}{2}\int_\mathcal
{X}\int
_\mathcal{X}\beta^2\gamma^{-\vert \mathbb{B}(x_1,R/2) \cup
\mathbb{B} (x_2,R/2) \vert}\mathbh{1}
\bigl\{\llVert x_1-x_2 \rrVert\leq R\bigr\}
\,dx_1 \,dx_2
\\
&&\qquad =c_{\Xi} \frac{e^{-\vert \mathcal{X} \vert}}{2} \int_\mathcal
{X}\int
_{\mathbb{B}(x_2,R) \cap\mathcal{X}}\beta^2 \gamma^{-\vert
\mathbb{B} (x_1,R/2) \vert}
\\
&&\hspace*{133pt}{}\times\gamma^{-\vert \mathbb{B}(x_2,R/2) \vert} \gamma^{\vert \mathbb
{B} (x_1,R/2)\cap\mathbb{B}(x_2,R/2) \vert} \,dx_1
\,dx_2
\\
&&\qquad \geq c_{\Xi} \frac{e^{-\vert \mathcal{X} \vert}}{2}\beta_0^2
\bigl\vert\mathcal{X}^{(-R)} \bigr\vert\int_{\mathbb{B}(0,R)}
\gamma^{\vert \mathbb{B}(x,R/2)\cap\mathbb{B}(0,R/2) \vert} \,dx,
\end{eqnarray*}
where we used that $\mathbb{B}$ denotes a ball in $\mathbb{R}^D$ and the
translation invariance of the Lebesgue measure. Since $\gamma\leq1$,
we have $u_\Xi(\xi)\le c_\Xi\beta_0^{\vert \xi\vert}$ for
all $\xi\in
\mathfrak{N}$. Integrating with respect to $\Po_1$ yields $c_\Xi
\ge\exp(\vert \mathcal{X} \vert(1-\beta_0))$. Thus one may choose
$\kappa= e^{-\beta_0 \vert \mathcal{X} \vert}\beta_0^2\vert
\mathcal{X}^{(-R_0)} \vert/2$.
\end{pf}

\subsection{Discrete processes}

Let $(\mathcal{X},\bolds{\alpha})$ be a general space with a diffuse measure
$\bolds{\alpha}$. Our aim is to compare Gibbs processes which live on a
finite subset $\Lambda=\Lambda_n=\{y_i\}_{i=1}^n$ of $\mathcal{X}$
with Gibbs processes on~$\mathcal{X}$. Let $V=\{V_i\}_{i=1}^n$ be a
partition of $\mathcal{X}$ such that $y_i\in V_i$ for all $i=1,\ldots,n$. A natural choice is the Voronoi tesselation, provided that
$\bolds{\alpha}(\partial V_i)=0$ for all $i=1,\ldots,n$. Define
$r_V=\max_{i=1,\ldots,n}\sup_{x\in V_i}\,d(x,y_i)$, the maximal radius
of the
cells in $V$. Furthermore, let $(\mathfrak{N}_\Lambda, \mathcal
{N}_\Lambda)$ denote the space of point measures on $\Lambda$ with
its natural $\sigma$-algebra, which coincides with the power set of
$\mathfrak{N}$.

Define the reference measure $\bolds{\alpha}_\Lambda$ on $\Lambda$ by
$\bolds{\alpha}_\Lambda(\{y_i\})=\bolds{\alpha}(V_i)$ for all
$i=1,\ldots,n$ and
let $\Po_\Lambda$ denote the Poisson process distribution on $\Lambda$
with intensity measure~$\bolds{\alpha}_\Lambda$. The Gibbs point processes
on $\Lambda$ are then defined in the obvious way, that is, as the
point processes that have a hereditary density with respect to
$\Po_\Lambda$.

Let $\Xi_\Lambda\sim$ $\gibbs(u_\Lambda)$ be a Gibbs process on
$\Lambda$. Define a point process $\Xi_U$ on $\mathcal{X}$ in the
following manner. Each point of $\Xi_\Lambda$ is replaced by a
\mbox{$\bolds{\alpha}(\cdot)|_{V_i}/\bolds{\alpha}(V_i)$-}distribu\-ted
point in the
corresponding cell $V_i$. More formally, if
$\Xi_\Lambda=\sum
_{i=1}^n N_i \delta_{y_i}$, then
%
%
\begin{equation}
\label{eqXXiU} \Xi_U=\sum_{i=1}^n
\sum_{l=1}^{N_i} \delta_{U_{il}},
\end{equation}
where the $U_{il}$ are all independent and $U_{il}\sim\bolds{\alpha
}(\cdot
)|_{V_i}/\bolds{\alpha}(V_i)$.

Define a function $t\colon\mathcal{X}\to\Lambda$ which maps each
point in $\mathcal{X}$ to its lattice point in $\Lambda$, that is,
$t(x) = y_i$ if $x \in V_i$ for $i=1,\ldots,n$. In the same spirit set
$t(\xi)=t ( \sum_{x\in\xi}\delta_{x} )=\sum_{x\in\xi
}\delta_{t(x)}$ for every $\xi\in\mathfrak{N}$.

%
\begin{lemma}
\label{lemdiscret-dens}
Let $\Xi_\Lambda\sim$ $\gibbs(u_\Lambda)$. Then the corresponding
point process $\Xi_U$ on $\mathcal{X}$ has density
$u_U(\xi)=u_\Lambda(t(\xi))$ with respect to $\Po_1$.
\end{lemma}

\begin{pf}
For $j_1,\ldots,j_n \in\mathbb{Z}_+$, set $k=\sum_{i=1}^n j_i$ and
\[
A_{j_1,\ldots,j_n}=\bigl\{\xi\in\mathfrak{N}\colon\xi(V_1)=j_1,
\ldots,\xi(V_n)=j_n\bigr\}.
\]

Then, writing
\[
\pmatrix{k\cr j_1,\ldots,j_n}=\frac{k!}{j_1!\cdots j_n!}
\]
for the multinomial coefficient, we have
%
%
\begin{eqnarray}
\label{eqPAi} \qquad&&{\mathbb P}(\Xi_U \in A_{j_1,\ldots,j_n})\nonumber
\\
&&\qquad ={\mathbb P}(\Xi_\Lambda\in A_{j_1,\ldots,j_n})\nonumber
\\
&&\qquad =\frac{e^{-\bolds{\alpha}(\mathcal{X})}}{k!}u_\Lambda\Biggl(\sum_{i=1}^n
j_i \delta_{y_i} \Biggr)\pmatrix{k\cr j_1,
\ldots,j_n} \bolds{\alpha} (V_1)^{j_1} \cdots
\bolds{\alpha}(V_n)^{j_n}\nonumber
\\
&&\qquad = \frac{e^{-\bolds{\alpha}(\mathcal{X})}}{k!}
\int_\mathcal
{X}\cdots\int
_\mathcal{X} u_\Lambda\Biggl(\sum
_{r=1}^k \delta_{t(x_r)}\Biggr)
\\
&&\hspace*{107pt}{}\times
\mathbh{1} \bigl\{ \#\{r\colon x_r\in V_1\}
=j_1,\ldots,\nonumber
\\
&&\hspace*{152pt}{}
 \#\{r\colon x_r\in V_n
\}=j_n \bigr\} \bolds{\alpha}(dx_1)\cdots\bolds{\alpha}
(dx_k)\nonumber
\\
&&\qquad ={\mathbb P}(\tXXi\in A_{j_1,\ldots,j_n}),
\nonumber
\end{eqnarray}
where $\tXXi\sim\gibbs( u_\Lambda\circ
t )$.

Note that the density $u_\Lambda(t(\cdot))$ is constant on any
$A_{j_1,\ldots,j_n}$. Hence, given that $\tXXi(V_1)=j_1,\ldots,\tXXi
(V_n)=j_n$, we may write $\tXXi$ as
$\sum_{i=1}^n\sum_{l=1}^{j_i}\delta_{U_{il}}$, where
the $U_{il}$ are all independent and $U_{il}\sim\bolds{\alpha}(\cdot
)|_{V_i}/\bolds{\alpha}(V_i)$.
Thus, for every measurable $h\colon\mathfrak{N}\to\mathbb{R}_+$ we
get ${\mathbb E}(h(\Xi_U)\mvert A_{j_1,\ldots,j_n})={\mathbb
E}(h(\tXXi)\mvert
A_{j_1,\ldots,j_n})$, which together with equation~(\ref{eqPAi}) and
the formula of total expectation yields the claim.
\end{pf}

Many Gibbs processes $\Xi$ on $\mathcal{X}$ with density $u$ have a
discrete analogon, which is obtained by restricting the density to
$\Lambda$ and renormalizing, that is, by using the unnormalized
density $\tilde{u}_\Lambda=u \vert_{\mathfrak{N}_\Lambda}$ on
$\mathfrak
{N}_\Lambda$, provided that it is $\bolds{\alpha}_\Lambda
$-integrable. Some
special care is required when evaluating $u$ at point configurations
$\xi\in\mathfrak{N}$ that have multi-points. For the continuous
Gibbs process such $\xi$ form a null set, whereas for the discrete
analogon the values of $\tilde{u}_\Lambda$ at such $\xi$ become important.
We avoid this problem by assuming that $u(\xi)=0$ for any $\xi$ with
multi-points, which leads to discrete analoga $\Xi_{\Lambda}$ that
may be represented as collections of Bernoulli random variables
$(I_y)_{y \in\Lambda}$.
Consequently $\Xi_U$ may have no more than one point in any cell~$V_i$.

By Theorem~\ref{thmmain2} and Lemma~\ref{lemdiscret-dens} we
immediately obtain the following proposition.
%
%
\begin{proposition}
\label{propdistU}
Suppose that $\Xi\sim\gibbs(\nu)$ satisfies $\SSS$. Let $\Xi
_\Lambda$ be its discrete analogon, and let $\Xi_U$ be given by
(\ref{eqXXiU}). Then
%
%
\begin{equation}
\quad d_{\mathrm{TV}}\bigl(\mathscr{L}(\Xi_U),\mathscr{L}(\Xi)\bigr)
\le c_1(\lambda)\int_\mathcal{X} {\mathbb E} \bigl| \nu
\bigl(t(x)\mvert t(\Xi_U)\bigr)-\nu(x\mvert\Xi_U) \bigr|
\bolds{\alpha}(dx).
\end{equation}
\end{proposition}

Consider the special case where $\Xi\sim\pip(\beta,\varphi)$ with
a constant $\beta$ and with $\varphi\leq1$ (inhibitory case). Our
process $\Xi_{\Lambda}$ is then an \emph{auto-logistic process} in
the terminology of \citet{Besag74}. In \citet{bmz82} convergence of the
corresponding $\Xi_U$-process density toward the $\Xi$-process
density is studied under a continuity condition on $\varphi$, without
providing rates.

We have
\[
\nu\bigl(t(x) \mvert t(\Xi_U)\bigr) = \beta\prod
_{y \in\Xi_U} \varphi\bigl(t(x),t(y)\bigr)\qquad\mbox{a.s.}
\]
Imitating the proof of Theorem~\ref{cordtv-pip}, we obtain from
Proposition~\ref{propdistU} the following result.

%
\begin{proposition}
\label{thmpipdiscrete}
Let $\Xi\sim\pip(\beta,\varphi)$ with constant $\beta$ and
$\varphi\leq1$.
Then
%
%
\begin{eqnarray}
\label{eqthmdiscr} && d_{\mathrm{TV}}\bigl(\mathscr{L}(\Xi_U),
\mathscr{L}(\Xi)\bigr)
\nonumber\\[-8pt]\\[-8pt]
&&\qquad \le c_1(\lambda) {\mathbb E}|\Xi_U| \beta\sup
_{y\in\mathcal
{X}}\int_\mathcal{X} \bigl| \varphi
\bigl(t(x),t(y)\bigr)-\varphi(x,y) \bigr| \bolds{\alpha}(dx).
\nonumber
\end{eqnarray}
\end{proposition}

%
\begin{corollary}
\label{corLcont}
Let $\Xi\sim\pip(\beta, \varphi)$ with constant $\beta$ and a
$\varphi\leq1$, that is, Lipschitz continuous with constant $L$ in
both components. Then
%
%
\begin{equation}
\label{eqcorL} d_{\mathrm{TV}}\bigl(\mathscr{L}(\Xi_U),
\mathscr{L}(\Xi)\bigr)\le2c_1(\lambda) {\mathbb E} |
\Xi_U| \beta L \bolds{\alpha}(\mathcal{X}) r_V.
\end{equation}
\end{corollary}
\begin{pf}
Note that by the triangle inequality,
\begin{eqnarray*}
&& \bigl|\varphi\bigl(t(x),t(y)\bigr)-\varphi(x,y)\bigr|
\\
&&\qquad \le \bigl|\varphi\bigl
(t(x),t(y)\bigr)-
\varphi\bigl(t(x),y\bigr)\bigr| + \bigl|\varphi\bigl(t(x),y\bigr)-\varphi(x,y)\bigr|
\le 2Lr_V.
\end{eqnarray*}\upqed
\end{pf}

%
\begin{example}
\label{exdisc-Strauss}
Let $\Xi$ be a Strauss process, that is, $\varphi(x,y) = \gamma
^{\mathbh{1}\{d(x,y)\le R\}}$ for a $\gamma\in[0,1]$. Write
$A(y,R_1,R_2)=\{x\in\mathcal{X}\colon R_1 < d(x,y)\le R_2\}$. For
$x\notin A(y,R-2r_V,R+2r_V)$, we have
\[
\gamma^{\mathbh{1}\{d(t(x),t(y))\le R\}}-\gamma^{\mathbh{1}\{
d(x,y)\le R\}}=0
\]
and for $x\in A(y,R-2r_V,R+2r_v)$ the modulus of the above difference
is at most $(1-\gamma)$. Hence by Proposition~\ref{thmpipdiscrete},
%
%
\begin{eqnarray}\label{eqexStrauss}
&& d_{\mathrm{TV}}\bigl(\mathscr{L}(\Xi_U),
\mathscr{L}(\Xi)\bigr)
\nonumber\\[-8pt]\\[-8pt]
&&\qquad \le c_1(\lambda) {\mathbb E}|\Xi
_U| \beta(1-\gamma)\sup_{y\in\mathcal{X}}\bolds{\alpha}
\bigl(A(y,R-2r_V,R+2r_V)\bigr).\nonumber
\end{eqnarray}
In the Euclidean case we get a linear rate in $r_V$, since
$\bolds{\alpha}(A(y,R-2r_V,R+2r_V))\le4\alpha_DD(R+2r_V)^{D-1}r_V$.
\end{example}

%
\begin{remark}
By combining the techniques of Corollary~\ref{corLcont} and
Example~\ref{exdisc-Strauss}, we obtain linear rates in $r_V$ for any
inhibitory interaction function $\varphi$ that is piecewise Lipschitz
continuous.
\end{remark}

To compute the distance between $\Xi$ and its discrete analogon $\Xi
_\Lambda$ we need another distance than the total variation. This is
because ${\mathbb P}(\Xi\subset\Lambda)={\mathbb P}(|\Xi|=0)$ and
thus $d_{\mathrm{TV}}(\mathscr{L}(\Xi_\Lambda),\mathscr{L}(\Xi))
\ge
1-{\mathbb P}(|\Xi|=0)$, whereas one would like to have a distance
that vanishes as $r_V \to0$. We use the following Wasserstein metric;
see \citet{bb92}, Section~3, for details. Let $\xi=\sum
_{i=1}^n\delta_{x_i}$ and $\eta=\sum_{i=1}^m\delta_{y_i}$. Define a
metric $d_1$ on $\mathfrak{N}$ by
\[
d_1(\xi,\eta)= %
\cases{ 1, &\quad if $n\neq m$,
\vspace*{4pt}\cr
\displaystyle\frac{1}{n}\min_{\sigma\in S_n}\sum
_{i=1}^n \min\bigl(d(x_i,y_{\sigma(i)}),1
\bigr), &\quad if $n=m$,}
\]
where $S_n$ denotes the permutation group of order $n$.
Denote by $\mathcal{F}_{2}$ the set of functions $f\colon\mathfrak
{N}\to[0,1]$ such that $|f(\xi)-f(\eta)| \le d_1(\xi,\eta)$ for
all $\xi,\eta\in\mathfrak{N}$. Our Wasserstein distance is then
defined by
%
%
\begin{equation}
\label{eqdef-Wasserstein} d_{2}\bigl(\mathscr{L}(\Xi),\mathscr
{L}(\mathrm{H})
\bigr)=\sup_{f\in
\mathcal
{F}_{2}}\bigl|{\mathbb E}f(\Xi)-{\mathbb E}f(\mathrm{H})\bigr|.
\end{equation}

We obtain the following theorem.
%
%
\begin{theorem}
\label{thmdiscrete-wasserstein}
Suppose that $\Xi\sim\gibbs(\nu)$ satisfies $\SSS$. Let $\Xi
_\Lambda$ be the discrete analogon. Then
%
%
\begin{equation}
\qquad d_{2}\bigl(\mathscr{L}(\Xi_\Lambda),\mathscr{L}(\Xi)\bigr)
\le r_V + c_1(\lambda) \int_\mathcal{X}
\mathbb{E} \bigl| \nu\bigl(t(x)\mvert t(\Xi_U)\bigr) - \nu(x\mvert
\Xi_U) \bigr| \bolds{\alpha}(dx).
\end{equation}
\end{theorem}
\begin{pf}
We have
\[
d_2\bigl(\mathscr{L}(\Xi_\Lambda),\mathscr{L}(\Xi)\bigr) \le
d_2\bigl(\mathscr{L}(\Xi_\Lambda),\mathscr{L}(
\Xi_U)\bigr) + d_2\bigl(\mathscr{L}(\Xi
_U),\mathscr{L}(\Xi)\bigr).
\]

For the first summand we obtain by the Lipschitz continuity of $f \in
\mathcal{F}_{2}$ that
\[
\sup_{f\in\mathcal{F}_{2}}\bigl|{\mathbb E}f(\Xi_\Lambda)-{\mathbb E}f(
\Xi_U)\bigr| \le{\mathbb E} \,d_1(\Xi_\Lambda,
\Xi_U)\le r_V,
\]
where we used that the distance between any point in $\Xi_\Lambda$
and its replacement point in $\Xi_U$ is at most $r_V$.

Since $\mathcal{F}_{2}\subset\mathcal{F}_{\mathrm{TV}}$, the
Wasserstein distance is
always bounded by the total variation distance. The second summand
above may therefore be bounded according to Proposition~\ref
{propdistU}, which yields the claim.
\end{pf}

\section{Couplings of spatial birth--death processes}
\label{seccoupling}

Let $b(\cdot\mvert\cdot),d(\cdot\mvert\cdot) \colon$ $\mathcal
{X}\times\mathfrak{N}\to\mathbb{R}_{+}$ be measurable functions
such that $\bar{b}(\xi):= \int b(x \mvert\xi) \bolds{\alpha
}(d x)$
$<\infty$ for every $\xi\in\mathfrak{N}$, and set $\bar{d}(\xi) =
\sum_{x \in\xi} \,d(x \mvert\xi)$ and $\bar{a}(\xi) = \bar{b}(\xi
) + \bar{d}(\xi)$. A \emph{spatial birth--death process} (SBDP)
$(Z(t))_{t\ge0}$ with \emph{birth rate} $b$ and \emph{death rate}
$d$ is a pure-jump Markov process on $\mathfrak{N}$ that can be
described as follows: given it is in state $\xi\in\mathfrak{N}$ and
$\bar{a}(\xi)>0$ it stays there for an $\Exp(\bar{a}(\xi
))$-distributed time, after which a point is added to $\xi$
(``birth'') with probability $\bar{b}(\xi) / \bar{a}(\xi)$ or deleted
from $\xi$ with probability $\bar{d}(\xi) / \bar{a}(\xi)$ (``death'').
If a birth occurs, the new point is positioned according to the density
$b(x \mvert\xi)/\bar{b}(\xi)$. If a death occurs, the point $x \in
\xi$ is omitted with probability $d(x \mvert\xi) / \bar{d}(\xi)$. In
the case $\bar{a}(\xi) = 0$ the SBDP is absorbed in $\xi$, that is,
stays there indefinitely. \citet{preston75} and \citet
{moellerwaage04}, Chapter~11 and Appendix~G, give more formal
definitions of general SBDP and a wealth of other results, including
conditions to assure that the SBDP is \emph{nonexplosive}, that is,
that with probability 1 only finitely many jumps can occur in any
bounded time interval. Denote by $(Z_{\xi}(t))_{t\ge0}$ the process
with deterministic starting configuration $\xi$, that is, $Z_{\xi}(0)
= \xi$.

We concentrate here on the case, where $b(x \mvert\xi):= \lambda(x
\mvert\xi)$ and $d(x \mvert\xi) = 1$ (``unit per-capita death
rate'') and where condition~$\SSS$ holds. Thus we obtain a nonexplosive
SBDP that is time-reversible with respect to $\gibbs(\lambda)$ and
converges in distribution to $\gibbs(\lambda)$ [see \citet
{moellerwaage04}, Propositions~\mbox{G.2--G.4}]. Time-reversibility with
respect to $\gibbs(\lambda)$ means that if $Z(0) \sim\gibbs(\lambda
)$, we have that $(Z(t))_{t \in[0,T]}$ and $(Z(T-t))_{t \in[0,T]}$
have the same distribution for any $T > 0$.
Since condition~$\SSS$ holds, we may also characterize the SBDP with
birth rate $\lambda$ and unit per-capita death rate as the unique
Markov process with infinitesimal generator
%
%
\begin{eqnarray}
\label{eqgen-Markov} \mathcal{A}h(\xi)&=&\int_{\mathcal{X}} \bigl[
h(\xi+
\delta_x)-h(\xi) \bigr] \lambda(x \mvert\xi) \bolds{\alpha}(d x)
\nonumber\\[-8pt]\\[-8pt]
&&{} +
\int_{\mathcal{X}} \bigl[ h(\xi-\delta_x)-h(\xi) \bigr]
\xi(d x)\nonumber
\end{eqnarray}
for all bounded measurable $h \colon\mathcal{X}\to\mathbb{R}$; see
\citet{ethierkurtz86}, Sections~4.2~and~4.11, Problem~5. The
domain $\mathscr{D}(\mathcal{A})$ of $\mathcal{A}$ is the set of all
measurable functions $h$ for which the right-hand side above is
well definied. $\mathscr{D}(\mathcal{A})$ contains at least all the
functions $h$ for which $\sup_{\xi\in\mathfrak{N}, x \in\mathcal
{X}} \vert h(\xi+\delta_x) - h(\xi) \vert$ is finite.

In what follows we construct a coupling $(Z_\xi(t),Z_\eta(t))_{t\ge
0}$ of two SBDPs with identical birth rate $\lambda$ started at
individual configurations $\xi,\eta\in\mathfrak{N}$. We introduce
the notation
\begin{eqnarray*}
\lambda_{\max} (x \mvert\xi,\eta) &=&\max\bigl(\lambda(x \mvert
\xi),
\lambda(x \mvert\eta)\bigr),\qquad \bar{\lambda}_{\max}(\xi,\eta)
=\int
_{\mathcal{X}} \lambda_{\max}(x\mvert\xi,\eta) \bolds{
\alpha}(dx),
\\
\lambda_{\min} (x \mvert\xi,\eta) &=&\min\bigl(\lambda(x \mvert
\xi),
\lambda(x \mvert\eta)\bigr),\qquad \bar{\lambda}_{\min}(\xi,\eta)
=\int
_{\mathcal{X}} \lambda_{\min
}(x\mvert\xi,\eta) \bolds{
\alpha}(dx).
\end{eqnarray*}

Define $(Z_\xi,Z_\eta)$ as a pure-jump Markov process with
right-continuous paths, holding intervals $D_1,D_2,\ldots,$ start time
$T_0:=0$, and jump times $T_j:=\sum_{i=1}^jD_i$ for all $j\ge1$.
Given $Z_\xi(T_{j-1})=\xi',Z_\eta(T_{j-1})=\eta'$ the distribution
of the next jump is described by the following random variables, which
are assumed to be independent of one another unless specified
otherwise. Let
\begin{eqnarray*}
D_j &\sim&\Exp\bigl(\bar{\lambda}_{\max}\bigl(
\xi',\eta'\bigr)+\bigl|\xi'\cup\eta
'\bigr|\bigr),
\\
G_j &\sim&\operatorname{Bernoulli} \biggl(\frac{\bar{\lambda
}_{\max}(\xi
',\eta
')}{\bar{\lambda}_{\max}(\xi',\eta')+|\xi'\cup\eta'|}
\biggr),
\\
Y_j &\sim&\frac{\lambda_{\max}(\cdot\mvert\xi',\eta')}{\bar
{\lambda}
_{\max}(\xi',\eta')},
\\
U_j &\sim&\operatorname{Unif}\bigl(\xi'\cup
\eta'\bigr),
\\
B_{\xi,j} \mid Y_j&\sim&\operatorname{Bernoulli} \biggl(
\frac{\lambda
(Y_j\mvert\xi')}{\lambda_{\max}(Y_j \mvert\xi',\eta')} \biggr),
\\
B_{\eta,j} \mid Y_j&\sim&\operatorname{Bernoulli} \biggl(
\frac{\lambda
(Y_j\mvert\eta')}{\lambda_{\max}(Y_j \mvert\xi',\eta')} \biggr),
\end{eqnarray*}
where $B_{\xi,j}$ and $B_{\eta,j}$ are maximally coupled given $Y_j$,
that is,
\[
\mathbb{P}(B_{\xi,j}=B_{\eta,j}=1 \mvert Y_j) =
\frac{\lambda_{\min}(Y_j \mvert\xi',\eta')}{\lambda_{\max
}(Y_j \mvert\xi',\eta')}.
\]
If $G_j=1$, set $Z_\xi(T_j)=Z_\xi(T_{j-1})+B_{\xi,j}\delta_{Y_j}$ and
$Z_\eta(T_j)=Z_\eta(T_{j-1})+B_{\eta,j}\delta_{Y_j}$.
If $G_j=0$, set $Z_\xi(T_j)=Z_\xi(T_{j-1})-\mathbh{1}\{U_j \in Z_\xi
(T_{j-1})\} \delta_{U_j}$ and $Z_\eta(T_j)=\break Z_\eta
(T_{j-1})- \mathbh{1}\{U_j \in Z_\eta(T_{j-1})\} \delta_{U_j}$.

In the special case where $\lambda(x\mvert\xi)$ does not depend on
the configuration $\xi$, that is, the time-reversible distribution is
the Poisson distribution with intensity function~$\lambda$, our
construction reduces to the coupling used in \citet{bb92}.

%
\begin{proposition}
\label{lemmacoupling}
Both components $Z_\xi,Z_\eta$ of the coupling are SBDPs with
generator (\ref{eqgen-Markov}).
\end{proposition}
\begin{pf}
Since by condition $\SSS$ the rate $a(\xi,\eta)=\bar{\lambda}_{\max
}(\xi,\eta)+|\xi\cup\eta|$ is bounded by $c+|\xi\cup\eta|$, for
a constant $c>0$, we get
\begin{eqnarray*}
{\mathbb P}(D_1 >t) &=&1-a(\xi,\eta)t+O\bigl(t^2\bigr),
\\
{\mathbb P}(D_1 \le t, D_1+D_2>t)&=&a(\xi,
\eta)t +O\bigl(t^2\bigr)
\end{eqnarray*}
and
\[
{\mathbb P}(D_1+D_2\le t)=O\bigl(t^2\bigr)
\]
as $t \to0$. Thus for a bounded function $h$
\begin{eqnarray*}
&& {\mathbb E} h\bigl(Z_\xi(t)\bigr)
\\[2pt]
&&\qquad = t\bigl(\bar{\lambda}_{\max}(\xi,\eta)+|\xi\cup\eta|\bigr)
\\[2pt]
&&\quad\qquad{} \times
\bigl({\mathbb E}h(\xi+B_{\xi,1}\delta_{Y_1})
{\mathbb P}(G_1=1) +{\mathbb E}h\bigl(\xi-\mathbh{1}\{U_1
\in\xi\}\delta_{U_1}\bigr){\mathbb P}(G_1=0) \bigr)
\\[2pt]
&&\quad\qquad{} + \bigl(1-t \bigl(\bar{\lambda}_{\max}(\xi,\eta)+|\xi
\cup\eta|\bigr) \bigr)h(\xi)+O\bigl(t^2\bigr).
\end{eqnarray*}
Since
\[
{\mathbb P}(G_1=1)=1-{\mathbb P}(G_1=0)=
\frac{\bar{\lambda}_{\max
}(\xi,\eta)}{\bar{\lambda}_{\max}(\xi,\eta)+|\xi\cup\eta|},
\]
we obtain
\begin{eqnarray*}
&& \lim_{t\to0} \frac{{\mathbb E}h(Z_\xi(t))-h(\xi)}{t}
\\[2pt]
&&\qquad =\bar{\lambda}_{\max}(\xi,\eta){\mathbb E} \bigl(h(
\xi+B_{\xi,1}\delta_{Y_1})-h(\xi) \bigr)
\\[2pt]
&&\quad\qquad{} +|\xi\cup\eta| {\mathbb E} \bigl(h\bigl(\xi-\mathbh{1}\{
U_1\in\xi\} \delta_{U_1}\bigr)-h(\xi) \bigr)
\\[2pt]
&&\qquad =\bar{\lambda}_{\max}(\xi,\eta)\int_{\mathcal{X}} \bigl[
h(\xi+\delta_y)-h(\xi) \bigr] \frac{\lambda(y \mvert\xi
)}{\lambda
_{\max}(y \mvert\xi, \eta)}\frac{\lambda_{\max}(y \mvert\xi,\eta)}{\bar{\lambda}_{\max}(\xi,\eta)}
\bolds{\alpha}(dy)
\\[2pt]
&&\quad\qquad{}+|\xi\cup\eta|\int_{\mathcal{X}} \bigl[ h\bigl(\xi
-\mathbh{1}\{u\in\xi\}\delta_u\bigr)-h(\xi) \bigr]
\frac{ (\xi
\cup\eta)(du)}{|\xi\cup\eta|}
\\[2pt]
&&\qquad =\int_{\mathcal{X}} \bigl[ h(\xi+\delta_y)-h(\xi)
\bigr] \lambda(y \mvert\xi) \bolds{\alpha}(dy)+\int_{\mathcal{X}}
\bigl[ h(\xi-\delta_u)-h(\xi) \bigr] \xi(du)
\\[2pt]
&&\qquad =\mathcal{A}h(\xi).
\end{eqnarray*}\upqed
\end{pf}

Define the coupling time as $\tau=\tau_{\xi,\eta}=\inf\{t \ge0: Z_\xi(t)=Z_\eta(t)\}$. In order to investigate the coupling time,
it is convenient to use the stopping times $\tau_0=0$, $\tau_k=\inf\{
t > \tau_{k-1}: Z_\xi(t)-Z_\eta(t) \neq Z_\xi(\tau
_{k-1})-Z_\eta(\tau_{k-1})\}$. These times are the times when
something interesting happens, that is, one of the noncommon points of
$Z_\xi$ and $Z_\eta$ dies or there is a birth in just one of the processes.

Let us call the event that a noncommon point dies a ``good death'' and
the event that only one process has a birth a ``bad birth.'' Assume
that there are $n$ noncommon points in $\xi$ and $\eta$, that is,
$\llVert \xi-\eta\rrVert=n$, where $\llVert \cdot\rrVert$
denotes the total
variation norm for signed measures. Define the event $A_n=\{\llVert
Z_\xi(\tau_1)-Z_\eta(\tau_1) \rrVert=n-1\}$ and the filtration
$\mathcal
{F}_t=\sigma( (Z_\xi(s),Z_\eta(s)); s \leq t )$.
Note that by construction $(Z_\xi,Z_\eta)$ has the strong Markov
property; see, for example, \citet{kallenberg02}, Theorem~12.14. An
easy calculation then gives us the following probabilities:
\begin{eqnarray*}
&& \mathbb{P}\bigl(\mbox{``good death'' at time }
T_j \mvert\mathcal{F}_{T_{j-1}}\bigr)
\\
&&\qquad =\frac{\llVert Z_{\xi}(T_{j-1})-Z_{\eta}(T_{j-1})
\rrVert}{\bar{\lambda}_{\max}(Z_\xi(T_{j-1}),Z_\eta
(T_{j-1}))+|Z_\xi
(T_{j-1})\cup Z_\eta(T_{j-1})|},
\\
&& \mathbb{P}\bigl(\mbox{``bad birth'' at time
}T_j\mvert\mathcal{F}_{T_{j-1}}\bigr)
\\
&&\qquad =\frac{\bar{\lambda}_{\max}(Z_\xi(T_{j-1}),Z_\eta
(T_{j-1}))-\bar{\lambda}_{\min}(Z_\xi(T_{j-1}),Z_\eta
(T_{j-1}))}{\bar{\lambda}
_{\max}(Z_\xi(T_{j-1}),Z_\eta(T_{j-1}))+|Z_\xi(T_{j-1})\cup Z_\eta
(T_{j-1})|}.
\end{eqnarray*}

%
\begin{lemma}
\label{lemmapn}
The probability of the event $A_n$ is bounded from below as
%
%
\begin{equation}
\label{eqpn} \mathbb{P}(A_n)\ge\biggl(1+\frac{1}{n}\sup
_{\|\xi'-\eta
'\|=n}\int_\mathcal{X}\bigl|\lambda\bigl(x\mid
\xi'\bigr)-\lambda\bigl(x\mid\eta'\bigr)\bigr| \bolds{
\alpha}(dx) \biggr)^{-1} >0.
\end{equation}
\end{lemma}
\begin{pf}
We argue in terms of the discrete time Markov chains $(Z_{\xi
}(T_j))_{j \in\mathbb{Z}_{+}}$ and $(Z_{\eta}(T_j))_{j \in\mathbb
{Z}_{+}}$ and refer to them as the \emph{jump chains} of our SBDPs.
Define the $\mathbb{N}$-valued random variable~$J$ by $\tau_1(\omega
) = T_{J(\omega)}(\omega)$ as the index of the first interesting
jump. Then
%
\begin{eqnarray*}
\mathbb{P}(A_n)
&=& \sum_{j=1}^{\infty} \int
_{\mathfrak{N}^2} \mathbb{P} \bigl( A_n \mvert J=j,
Z_{\xi}(T_{j-1}) = \xi', Z_{\eta}(T_{j-1})
= \eta' \bigr)
\\
&&\hspace*{31pt}{}\times  \mathbb{P} \bigl( J=j, Z_{\xi}(T_{j-1}) \in d
\xi', Z_{\eta}(T_{j-1}) \in d \eta'
\bigr)
\\
&=& \sum_{j=1}^{\infty} \int
_{\mathfrak{N}^2} \mathbb{P} \bigl(\bigl\llVert Z_{\xi
}(T_j)-Z_{\eta}(T_j)
\bigr\rrVert= n-1 \big\mvert
\\
&&\hspace*{46pt} \bigl\llVert Z_{\xi}(T_i)-Z_{\eta}(T_i)
\bigr\rrVert= n,\ \forall i \leq j-1,
\\
&&\hspace*{46pt} \bigl\llVert Z_{\xi}(T_j)-Z_{\eta}(T_j)
\bigr\rrVert\in\{n-1,n+1\},
\\
&&\hspace*{76pt} Z_{\xi}(T_{j-1}) = \xi',
Z_{\eta}(T_{j-1}) = \eta' \bigr)
\\
&&\hspace*{31pt}{}\times \mathbb{P} \bigl(J=j, Z_{\xi}(T_{j-1}) \in d
\xi', Z_{\eta}(T_{j-1}) \in d \eta'
\bigr)
\\
&=& \sum_{j=1}^{\infty} \int
_{\mathfrak{N}^2} \frac{n}{n+(\bar
{\lambda}
_{\max}(\xi',\eta')-\bar{\lambda}_{\min}(\xi',\eta'))}
\\
&&\hspace*{32pt}{}\times  \mathbb{P} \bigl(J=j, Z_{\xi}(T_{j-1}) \in d
\xi', Z_{\eta}(T_{j-1}) \in d \eta'
\bigr)
\\
&=& \mathbb{E} \biggl( \frac{n}{n+(\bar{\lambda}_{\max}(Z_{\xi
}(T_{J-1}),Z_{\eta
}(T_{J-1}))-\bar{\lambda}_{\min}(Z_{\xi}(T_{J-1}),Z_{\eta
}(T_{J-1})))} \biggr)
\\
&\ge&\frac{n}{n+\sup_{\llVert \xi'-\eta' \rrVert=n}(\bar
{\lambda}_{\max}(\xi
',\eta')-\bar{\lambda}_{\min}(\xi',\eta'))}.
\end{eqnarray*}
The claim follows since
\[
\bar{\lambda}_{\max}(\xi,\eta)-\bar{\lambda}_{\min}(\xi,\eta
)=\int_\mathcal{X}\bigl|\lambda(x\mid\xi)-\lambda(x\mid\eta)\bigr| \bolds
{\alpha}(dx),
\]
which is uniformly bounded in $\xi$ and $\eta$ by condition $\SSS$.
\end{pf}

%
\begin{theorem}
\label{thmeps-rec}
For all configurations $\xi,\eta$ the coupling time $\tau_{\xi,\eta
}$ is integrable. In particular if $\xi$ and $\eta$ differ in only
one point, we have for any $n^* \in\mathbb{N}\cup\{\infty\}$,
%
%
\begin{eqnarray}\label{eqrec-eps-sol}
\mathbb{E}\tau_{\xi,\eta} &\le&\bigl
(n^*-1\bigr)! \biggl(
\frac{\varepsilon
}{c} \biggr)^{n^*-1} \Biggl(\frac{1}{c}\sum
_{i=n^*}^\infty\frac{c^i}{i!}+\int
_0^c\frac{1}{s}\sum
_{i=n^*}^\infty\frac{s^i}{i!} \,ds \Biggr)
\nonumber\\[-8pt]\\[-8pt]
&&{} +\frac{1+\varepsilon}{\varepsilon}\sum_{i=1}^{n^*-1}
\frac{\varepsilon^{i}}{i},\nonumber
\end{eqnarray}
where
\[
\varepsilon =\sup_{\|\xi-\eta\|=1}\int_\mathcal{X}\bigl|
\lambda(x\mid\xi)-\lambda(x \mid\eta)\bigr| \bolds{\alpha}(dx)
\]
and
\[
c = c\bigl(n^*\bigr)=\sup_{\|\xi-\eta\|\ge n^*}\int_\mathcal{X}\bigl|
\lambda(x\mid\xi)-\lambda(x \mid\eta)\bigr| \bolds{\alpha}(dx)
\]
with the interpretations detailed in Theorem~\ref{thmmain2}.
The constants $\varepsilon$ and $c$ are finite by condition $\SSS$.
\end{theorem}

The following lemma treats the case $n^*=1$ and will be useful for the
proof of Theorem~\ref{thmeps-rec}.

%
\begin{lemma}
\label{lemmaTint} For all configurations $\xi,\eta$ the coupling
time $\tau_{\xi,\eta}$ is integrable. In particular if $\xi$ and
$\eta$ differ in only one point, we have
%
%
\begin{equation}
\mathbb{E}\tau_{\xi,\eta} \le\frac{e^c-1}{c} +\int_0^c
\frac{e^s-1}{s} \,ds,
\end{equation}
where $c=\sup_{\xi,\eta\in\mathfrak{N}} \int\vert \lambda(x
\vert\xi)-\lambda(x \vert\eta) \vert \bolds{\alpha}(dx)$,
which is finite by
$\SSS$.
\end{lemma}

\begin{pf}
Let $p_n=(1+c/n)^{-1}$, $n\ge1$. Construct a new pure-jump Markov
process $(Y(t))_{t\ge0}$ on $\mathbb{Z}_+$ by the following rule.
Given $Y(t)$ is in state $n \in\mathbb{Z}_{+}$, after an
exponentially distributed time with mean $1/n$, it jumps to $n-1$ with
probability $p_n$ and to $n+1$ with probability $1-p_n$. Define
stopping times $\tilde{\tau}_n=\inf\{t\ge0: Y_n(t)=0 \}$ for
all $n\ge0$, where $(Y_n(t))_{t \ge0}$ denotes the process started at $n$.

We show that $\tau_{\xi,\eta}$ is stochastically dominated by
$\tilde{\tau}_n$, and therefore $\mathbb{E}\tau_{\xi,\eta} \leq
\mathbb{E}
\tilde{\tau}_n$. Denote by $X=\|Z_\xi-Z_\eta\|$ the process
counting the noncommon points, and define a new jump process\vadjust{\goodbreak} $Y'$ by
the following construction. Set $Y'(0)=X(0)$. If $X$ and $Y'$ are on
the same level, say $n$, they move together. If $X$ jumps, then $Y'$
jumps with probability $p_n$ to $n-1$ and with probability $1-p_n$ to
$n+1$. Since ${\mathbb P}(A_n) \ge p_n$, the jumps can be coupled such
that $Y'$ stays above $X$. If they are separated, let $Y'$ behave like
$Y$ until they meet again. Thus we have $Y'(t)\ge X(t)$ for all $t\ge
0$, and hence $\tau_{\xi,\eta}$ is stochastically dominated by
$\tilde{\tau}'_n = \inf\{t\ge0: Y'_n(t)=0 \}$, where
$(Y'_n(t))_{t \ge0}$ denotes the process $Y'$ started at $n$.

Note that $Y'$ has the same transition probabilities as $Y$, but its
holding times are sometimes those of $X$.
Since the processes $Z_\xi$ and $Z_\eta$ have unit per-capita death
rates, each of their noncommon points dies independently after a
standard exponentially distributed time. If there are $n$ such points,
the minimum of these times is exponentially distributed with mean
$1/n$. Hence the holding times of $X$ at $n$, and therefore also all of
the holding times of $Y'$ at $n$, may be coupled with exponentially
distributed random variables with mean $1/n$ that are almost surely
larger or equal. This yields a coupling of $Y'$ and $Y$ with exactly
matched jump chains where $Y$ is just a slower version of $Y'$. Hence
$\tilde{\tau}'_n$ is stochastically dominated by~$\tilde{\tau}_n$.

Define now $e_n={\mathbb E}\tilde{\tau}_n$, $n\ge0$. Then, by conditioning
on the next jump in the $Y$-chain, we obtain
%
%
\begin{eqnarray} \label{eqrec-eq}
e_n&=& \biggl( e_{n-1}+\frac{1}{n} \biggr)
p_n+ \biggl( e_{n+1}+\frac
{1}{n} \biggr)
(1-p_n)
\nonumber\\[-8pt]\\[-8pt]
&=& e_{n-1}p_n+e_{n+1}(1-p_n)+
\frac{1}{n}.\nonumber
\end{eqnarray}
If $c=0$, we have $e_n=e_{n-1}+1/n=\sum_{i=1}^n1/i$, in particular
$e_1=1$. For $c>0$, define $a_n=e_n-e_{n-1}$, $n\ge1$, assuming that
the $e_n$ are finite. Then $p_n=(1+c/n)^{-1}$ yields
%
%
\begin{equation}
\label{eqrec} a_{n+1}=\frac{n}{c}a_n-
\frac{1}{c}-\frac{1}{n} \qquad\mbox{for all } n\ge1.
\end{equation}
Since $e_0=0$, the starting point is $a_1=e_1$. The general solution of
(\ref{eqrec}) is given by
%
%
\begin{equation}
\label{eqrec-sol} a_n=\sum_{i=0}^\infty
\frac{c^i}{\prod_{k=0}^i(n+k)} \biggl(1+\frac
{c}{n+i} \biggr)+C\frac{(n-1)!}{c^{n-1}}
\end{equation}
for an arbitrary constant $C\in\mathbb{R}$.

A result in \citet{gs01} [Exercise 6, page 265] states that the
sequence of expected return times is the smallest nonnegative solution
of (\ref{eqrec-eq}), which yields that the $e_n$ are in fact finite
and given by setting $C=0$. We obtain for all $n\ge1$ that
%
%
\begin{equation}
\label{eqen} e_n =\sum_{i=1}^na_n
\le na_1=ne_1
\end{equation}
and
\[
e_1=a_1=\sum_{i=0}^\infty
\frac{c^i}{(1+i)!} \biggl(1+\frac
{c}{1+i} \biggr)=\frac{e^c-1}{c} +\int
_0^c\frac{e^s-1}{s} \,ds.
\]
Note that $e_1$ converges to $1$ for $c$ going to $0$.
\end{pf}

\begin{pf*}{Proof of Theorem~\ref{thmeps-rec}}
Let $\xi$ and $\eta$ be point configurations differing in $n$ points,
that is, they can be written as
\[
\xi=\zeta+\sum_{i=1}^k
\delta_{y_i} \quad\mbox{and}\quad\eta=\zeta+\sum
_{i=k+1}^n\delta_{y_i},
\]
where $0\le k\le n$, $y_1,\ldots,y_n$ are the noncommon points of
$\xi$ and $\eta$, and $\zeta=\xi\cap\eta$. Then
\begin{eqnarray*}
\bigl|\lambda(x \mvert\xi)-\lambda(x \mvert\eta)\bigr|&=& \Biggl| \sum
_{j=1}^k \Biggl[ \lambda\Biggl( x \bigg\mvert \zeta+
\sum_{i=1}^j\delta_{y_i}
\Biggr)- \lambda\Biggl(x \bigg\mvert \zeta+ \sum_{i=1}^{j-1}
\delta_{y_i} \Biggr) \Biggr]
\\
&&\hspace*{2pt}{} -\sum_{j=k+1}^n \Biggl[ \lambda
\Biggl(x \bigg\mvert \zeta+ \sum_{i=k+1}^j
\delta_{y_i} \Biggr)- \lambda\Biggl(x \bigg\mvert \zeta+ \sum
_{i=k+1}^{j-1}\delta_{y_i} \Biggr) \Biggr] \Biggr|.
\end{eqnarray*}
By the triangle inequality we obtain
\begin{eqnarray*}
&& \sup_{\|\xi-\eta\|=n}\int_\mathcal{X}\bigl|\lambda(x\mid\xi
)-\lambda(x\mid\eta)\bigr| \bolds{\alpha}(dx)
\\
&&\qquad \le n\sup_{\|\xi-\eta\|=1}\int_\mathcal{X}\bigl|\lambda(x
\mid\xi)-\lambda(x\mid\eta)\bigr| \bolds{\alpha}(dx)=n\varepsilon.
\end{eqnarray*}
Thus by Lemma~\ref{lemmapn} we have ${\mathbb P}(A_n) \ge
(1+\varepsilon)^{-1}$ for $n\ge1$ and ${\mathbb P}(A_n) \ge
(1+c/n)^{-1}$ for $n \ge n^*$. Assume $\varepsilon>0$, $c>0$ and
$n^{*} \in\mathbb{N}$. Replace the jump-down probabilities $p_n$ of
$Y$ in the proof of Lemma~\ref{lemmaTint} by the above bounds. We
then obtain for the differences $a_n=e_n-e_{n-1}$ of the expected
return times to zero the recursion equations
%
%
\begin{equation}
\label{eqrec-eps} %
\cases{ \displaystyle a_{n+1}=
\frac{1}{\varepsilon}a_n-\frac{1+\varepsilon
}{n\varepsilon}, &\quad for $1\le n < n^*$,
\vspace*{4pt}\cr
\displaystyle a_{n+1}=\frac{n}{c}a_n -
\frac{1}{c}-\frac{1}{n}, &\quad for $n\ge n^*$.}
\end{equation}
The differences for larger $n$ must still be the same. Hence the proof
of Lemma~\ref{lemmaTint} gives
\[
a_n=\sum_{i=0}^\infty
\frac{c^i}{\prod_{k=0}^i(n+k)} \biggl(1+\frac
{c}{n+i} \biggr)
\]
for all $n\ge n^*$. The second recursion equation in (\ref
{eqrec-eps}) is best solved backwards, which yields the solution
\[
a_{n^*-k}=\varepsilon^k a_{n^*} +(1+\varepsilon)
\sum_{i=1}^k\frac
{\varepsilon^{k-i}}{n^*-i}
\]
for all $1\le k \le n^*-1$.
Thus
%
%
\begin{equation}
\label{eqe1fordummies}
\qquad e_1=a_1=\varepsilon^{n^*-1}\sum
_{i=0}^\infty\frac{c^i}{\prod_{k=0}^i(n^*+k)} \biggl(1+
\frac{c}{n^*+i} \biggr)+ \frac{1+\varepsilon}{\varepsilon}\sum
_{i=1}^{n^*-1}
\frac
{\varepsilon^{n^*-i}}{n^*-i},
\end{equation}
which can be rewritten as (\ref{eqrec-eps-sol}).

Letting $n^{*} \to\infty$ in the inequality $\mathbb{E}\tau_{\xi,\eta}
\leq e_1$, we obtain $\mathbb{E}\tau_{\xi,\eta} \leq\frac
{1+\varepsilon}{\varepsilon
} \sum_{i=1}^\infty\frac{\varepsilon^i}{i}$ if $\varepsilon< 1$
irrespective of
$c$, which justifies setting the long first summand in (\ref
{eqrec-eps-sol}) to zero. Analogously as in the proof of Lemma~\ref
{lemmaTint} we get for $\varepsilon=0$ (which implies $c=0$) that
$e_1=1$, and for $\varepsilon>0, c=0$ that $a_{n^*}=1/n^*$, which
justifies the interpretation of the upper bound in (\ref
{eqrec-eps-sol}) in the limit sense.

If $c$ does not depend on $n^*$, we choose
$n^*$ such that $(1+\varepsilon)^{-1}>(1+c/n)^{-1}$ for $n<n^*$ and
$(1+\varepsilon)^{-1}\le(1+c/n)^{-1}$ for $n\ge n^*$. This is
obviously the optimal choice and leads to
$n^*=\lceil c/\varepsilon\rceil$.
\end{pf*}

One can also couple SBDPs with random starting configurations. It is
convenient to use notation of the form $Z_{\Xi}$ also if $\Xi$ is a
point process with the obvious meaning that $(Z_{\Xi}(t))_{t \geq0}$
is an SBDP with generator~(\ref{eqgen-Markov}) and $Z(0) = \Xi$
almost surely. Since this may lead to confusing notation when dealing
with two processes, we always distinguish the processes by adding a
prime, thus writing $(Z_{\Xi}(t))$ and $(Z'_{\tXXi}(t))$ for point
processes $\Xi$ and~$\tXXi$.

%
\begin{proposition}
\label{thmcouplingH}
Assume that $\Xi$ and $\tXXi$ are Gibbs processes satisfying $\SSS$.
Consider the coupling $(Z_\Xi(t),Z'_{\tXXi}(t))$, $t\ge0$. Then the
coupling time $\tau_{\Xi,\tXXi}=\inf\{t\ge0: Z_\Xi
(t)=Z'_{\tXXi}(t)\}$ is integrable.
\end{proposition}
\begin{pf}
The Georgii--Nguyen--Zessin equation and condition~$\SSS$ yield
\[
{\mathbb E}|\Xi|={\mathbb E}\int_\mathcal{X}1 \Xi(dx)= {\mathbb
E}\int_\mathcal{X}\nu(x\mvert\Xi) \bolds{\alpha}(dx)<\infty,
\]
where $\nu$ is the conditional intensity of $\Xi$, and analogously
${\mathbb E}|\tXXi| <\infty$. Then (\ref{eqen}) implies
\begin{eqnarray*}
{\mathbb E}\tau_{\Xi,\tXXi} &=&{\mathbb E}\bigl({\mathbb E}\bigl(
\tau_{\Xi,\tXXi} \mvert|\Xi|,|\tXXi|\bigr)\bigr)
\\
&\le& {\mathbb E}\bigl({\mathbb
E}\bigl(e_{|\Xi|+|\tXXi|}\mvert|\Xi|,|\tXXi|\bigr)\bigr) \le e_1{\mathbb E}\bigl(|
\Xi|+|\tXXi|\bigr) < \infty.
\end{eqnarray*}\upqed
\end{pf}

In particular the coupling time $\tau_{\xi,\mathrm{H}}$ for (nonrandom)
$\xi\in\mathfrak{N}$ and $\mathrm{H}$ is integrable.

\section{Stein's method for Gibbs process approximation}
\label{secSteinforGibbs}

Stein's method, originally conceived for normal approximation
[\citet
{stein72}], has evolved over the last forty years to become an
important tool in many areas of probability theory and for a wide range
of approximating distributions. See \citet{steinintro} for an overview.

A milestone in the evolution of Stein's method was the discovery in
\citet{barbour88} that a natural Stein equation may often be set up by
choosing as a right-hand side the infinitesimal generator of a Markov
process whose stationary distribution is the approximating distribution
of interest. Many important developments stem from this so-called \emph
{generator approach} to Stein's method, and several of them concern
point process approximation, such as \citet{bb92}, \citet{barman02},
\citet{schumi4} or \citet{xiazhang12}.

In this section we develop the generator approach for Gibbs process
approximation.
Let $\mathrm{H}\sim\gibbs(\lambda)$ be our approximating Gibbs processes
satisfying $\SSS$. Define the generator
%
%
\begin{equation}\label{eqgen}
\mathcal{A}h(\xi) =\int_{\mathcal{X}} \bigl[ h(\xi+
\delta_x)-h(\xi) \bigr] \lambda(x \mvert\xi) \bolds{\alpha}(d x)
+\int_{\mathcal{X}} \bigl[ h(\xi-\delta_x)-h(\xi) \bigr]
\xi(d x)\hspace*{-30pt}
\end{equation}
for all $h \colon\mathfrak{N}\to\mathbb{R}$ in its domain $\mathscr
{D}(\mathcal{A})$. In Section~\ref{seccoupling} we noted that
$\mathcal{A}$ is the generator of a spatial birth--death process $Z$
with stationary distribution $\gibbs(\lambda)$, and that its domain
contains at least all functions $h$ with bounded first differences.

For any measurable $f \colon\mathfrak{N}\to\mathbb{R}$ set up the
so-called \emph{Stein equation} formally as
%
%
\begin{equation}
\label{eqstein} f(\xi)-\mathbb{E}f(\mathrm{H})=\mathcal{A}h(\xi).
\end{equation}
A first goal is to find a function $h=h_f$ that satisfies this
equation. By analogy to the Poisson process case, a natural candidate
is given by
%
%
\begin{equation}
\label{eqstein-solution} h_f(\xi)=-\int_0^\infty
\bigl[ {\mathbb E}f\bigl(Z_\xi(t)\bigr)-{\mathbb E}f(\mathrm{H})
\bigr] \,dt.
\end{equation}

The following lemma shows that $h_f$ is indeed a solution to
equation~(\ref{eqstein}) if $f \in\mathcal{F}_{\mathrm{TV}}$.

%
\begin{lemma}
\label{lemmastein-sol-b}
Assume that $f$ is bounded and $\mathrm{H}$ satisfies $\SSS$. Then
$h_f$ is
well defined, that is, the integral exists for all $\xi\in\mathfrak
{N}$, and it is a solution to~(\ref{eqstein}).
\end{lemma}

\begin{pf} 
We use the coupling~$(Z_{\xi},Z'_{\mathrm{H}})$ from Section~\ref
{seccoupling}, where $Z'_{\mathrm{H}}$ is started in the random
configuration $\mathrm{H}$, as explained at the end of that section. Since
$\mathscr{L}(\mathrm{H})$ is the stationary measure\vadjust{\goodbreak} of the SBDP $Z$, we
have ${\mathbb E}f(\mathrm{H})={\mathbb E}f(Z'_{\mathrm{H}}(t))$ for
all $t\ge0$. Thus
\begin{eqnarray*}
\int_0^\infty\bigl|{\mathbb E}f\bigl(Z_\xi(t)
\bigr)-{\mathbb E}f(\mathrm{H})\bigr| \,dt
&=& \int_0^\infty{
\mathbb E}\bigl|f\bigl(Z_\xi(t)\bigr)-f\bigl(Z'_{\mathrm{H}}(t)
\bigr)\bigr| \mathbh{1} \{\tau_{\xi,\mathrm{H}} >t\} \,dt
\\
&\le& 2 \llVert f \rrVert_\infty\int_0^\infty{
\mathbb P}(\tau_{\xi,\mathrm{H}} >t) \,dt
\\
&=& 2 \llVert f \rrVert_\infty{\mathbb E}(\tau_{\xi,\mathrm{H}})
<\infty
\end{eqnarray*}
by Proposition~\ref{thmcouplingH}. Hence $h_f$ is well defined.
The Markov property of the SBDP implies $\mathscr{L} ( Z'_{Z_\xi
(s)}(t) ) = \mathscr{L} ( Z_\xi(t+s) )$. Thus by the
substitution $v=t+s$
\begin{eqnarray*}
&& \frac{1}{s} \bigl({\mathbb E}h_f\bigl(Z_\xi(s)
\bigr)-h_f(\xi) \bigr)
\\
&&\qquad =\frac{1}{s} \biggl(-\int
_0^\infty{\mathbb E} \bigl[f\bigl(Z_\xi(t+s)
\bigr)-f(\mathrm{H}) \bigr] \,d t
+\int_0^\infty{\mathbb E}
\bigl[f\bigl(Z_\xi(t)\bigr)-f(\mathrm{H}) \bigr] \,d t \biggr)
\\
&&\qquad = \frac{1}{s}\int_0^s{\mathbb E}
\bigl[f\bigl(Z_\xi(v)\bigr)-f(\mathrm{H}) \bigr] \,d v.
\end{eqnarray*}
By condition $\SSS$ we have ${\mathbb P}(Z_\xi(v)\neq\xi)=O(v)$ and
since $f$ is bounded, ${\mathbb E}f(Z_\xi(v))= f(\xi)+ O(v)$, which implies
\[
\frac{1}{s}\int_0^s{\mathbb E}f
\bigl(Z_\xi(v)\bigr) \,d v=f(\xi)+O(s).
\]
Thus
\[
\mathcal{A}h_f(\xi)=\lim_{s \to0}\frac{1}{s}
\bigl({\mathbb E}h_f\bigl(Z_\xi(s)
\bigr)-h_f(\xi) \bigr)=f(\xi)-{\mathbb E}f(\mathrm{H}).
\]\upqed
\end{pf}

Define then the \emph{Stein factor} as
%
%
\begin{equation}
\label{eqc1def} c_1(\lambda)=\sup_{f\in\mathcal{F}_{\mathrm
{TV}}}\,\sup
_{x\in
\mathcal{X}, \xi\in
\mathfrak{N}}\bigl|h_f(\xi+\delta_x)-h_f(
\xi)\bigr|.
\end{equation}

We are now ready to give proofs for the results in Section~\ref{ssecmain}.
\begin{pf*}{Proof of Theorem~\ref{thmmain2}}
By the Stein equation (\ref{eqstein})
\[
d_{\mathrm{TV}}\bigl(\mathscr{L}(\Xi),\mathscr{L}(\mathrm
{H})\bigr)=\sup
_{f\in\mathcal{F}_{\mathrm{TV}}} \bigl\vert{\mathbb E}f(\Xi
)-{\mathbb E}f(\mathrm{H})
\bigr\vert= \sup_{f\in\mathcal{F}_{\mathrm{TV}}} \bigl\vert
{\mathbb E}
\mathcal{A}h_f(\Xi) \bigr\vert,
\]
where $h_f$ is the Stein solution given in (\ref{eqstein-solution}).
Then the Georgii--Nguyen--Zessin equation~(\ref{eqgnz}) yields
%
%
\begin{eqnarray}
\label{eqgenbound}
&& \bigl\vert\mathbb{E}\mathcal{A} h_f(\Xi) \bigr\vert\nonumber
\\
&&\qquad = \biggl| \mathbb{E}\int_{\mathcal{X}} \bigl[ h_f(\Xi+
\delta_x) - h_f(\Xi) \bigr]\lambda(x \mvert\Xi)
\bolds{\alpha}(d x)\nonumber
\\
&&\hspace*{53pt}{}+ \mathbb{E}\int_{\mathcal{X}} \bigl[
h_f(\Xi- \delta_x) - h_f(\Xi) \bigr]
\Xi(d x) \biggr|
\nonumber\\[-8pt]\\[-8pt]
&&\qquad = \biggl| \mathbb{E}\int_{\mathcal{X}} \bigl[
h_f(\Xi+ \delta_x) - h_f(\Xi) \bigr]
\bigl(\lambda(x \mvert\Xi)-\nu(x\mvert\Xi)\bigr) \bolds{\alpha
}(d x) \biggr|\nonumber
\\
&&\qquad \le \sup_{\xi,\eta\in\mathfrak{N}, \llVert \xi-\eta\rrVert
=1}\bigl|h_f(
\xi)-h_f(\eta)\bigr|\int_\mathcal{X}\mathbb{E}\bigl|\nu(x
\mvert\Xi)-\lambda(x \mvert\Xi)\bigr| \bolds{\alpha}(dx)\nonumber
\\
&&\qquad \leq c_1(\lambda)\int_\mathcal{X}
\mathbb{E}\bigl|\nu(x\mvert\Xi)-\lambda(x \mvert\Xi)\bigr| \bolds{\alpha}(dx).\nonumber
\end{eqnarray}

Consider the coupling $(Z_{\xi+\delta_x}(t),Z_{\xi}(t))_{t\ge0}$ in
the sense of Section~\ref{seccoupling}. Then by equation~(\ref
{eqstein-solution})
\begin{eqnarray*}
\bigl|h_f(\xi+\delta_x)-h_f(\xi)\bigr|&=&\biggl
\vert{\mathbb E}\int_0^\infty\bigl[f
\bigl(Z_{\xi+\delta_x}(t)\bigr)-f\bigl(Z_\xi(t)\bigr) \bigr]
\mathbh{1}\{\tau_{\xi
+\delta
_x,\xi}>t\} \,dt\biggr\vert
\\
&\le&\sup_{\xi,\eta\in\mathfrak{N}}\bigl|f(\xi)-f(\eta)\bigr|\int
_0^\infty{
\mathbb P}(\tau_{\xi+\delta_x,\xi}>t) \,dt \le{\mathbb E}\tau_{\xi
+\delta_x,\xi},
\end{eqnarray*}
where we used that $0\le f\le1$. The upper bound on $c_1(\lambda)$
follows now from Theorem~\ref{thmeps-rec}.
\end{pf*}

%
\begin{remark}
\label{remarkesssup}
In the proof above the $\sup_{\xi,\eta\in\mathfrak{N}, \llVert
\xi-\eta\rrVert=1}$ could actually be replaced by an essential supremum
with respect to $\mathscr{L}(\Xi)+\mathscr{L}(\mathrm{H})$. This
can be
seen as follows. Let $N \in\mathcal{N}$ be a null set with respect to
both $\mathscr{L}(\Xi)$ and $\mathscr{L}(\mathrm{H})$. Without loss of
generality we may set the densities of $\Xi$ and $\mathrm{H}$ to zero on
$N$. By the hereditarity of the densities we have
\begin{eqnarray*}
&& \bigl[ h_f(\Xi+ \delta_x) - h_f(\Xi)
\bigr]\bigl(\lambda(x \mvert\Xi)-\nu(x\mvert\Xi)\bigr)
\\
&&\qquad = \bigl[ h_f(\Xi+ \delta_x) - h_f(
\Xi) \bigr] \mathbh{1}\bigl\{ \Xi+ \delta_x \in N^c, \Xi
\in N^c\bigr\} \bigl(\lambda(x \mvert\Xi)-\nu(x\mvert\Xi)\bigr),
\end{eqnarray*}
whence it follows that the first inequality of (\ref{eqgenbound})
holds also for the essential supremum.

A consequence of this replacement is that it suffices to take the
essential supremum with respect to~$\mathscr{L}(\Xi)+\mathscr
{L}(\mathrm{H})$ for the computation of the constants $\varepsilon$
and $c$
in Theorem~\ref{thmmain2}.
\end{remark}

\begin{pf*}{Proof of Theorem~\ref{cordtv-pip}}
Let $\nu$ and $\lambda$ denote the conditional intensities of $\Xi$
and $\mathrm{H}$. Then for $\xi= \sum_{i=1}^n \delta_{y_i} \in
\mathfrak
{N}$ we have
\begin{eqnarray*}
&& \nu(x \mvert\xi) - \lambda(x \mvert\xi)
\\
&&\qquad = \beta(x) \biggl( \prod_{y \in\xi}
\varphi_1(x,y) - \prod_{y \in\xi}
\varphi_2(x,y) \biggr)
\\
&&\qquad = \beta(x) \sum_{j=1}^n \Biggl( \Biggl(
\prod_{i=1}^j \varphi
_1(x,y_i) \Biggr) \Biggl( \prod
_{i=j+1}^n \varphi_2(x,y_i)
\Biggr)
\\
&&\hspace*{77pt}{} - \Biggl( \prod_{i=1}^{j-1}
\varphi_1(x,y_i) \Biggr) \Biggl( \prod
_{i=j}^n \varphi_2(x,y_i)
\Biggr) \Biggr)
\\
&&\qquad = \beta(x) \sum_{j=1}^n \Biggl( \bigl(
\varphi_1(x,y_j) - \varphi_2(x,y_j)
\bigr) \Biggl( \prod_{i=1}^{j-1}
\varphi_1(x,y_i) \Biggr) \Biggl( \prod
_{i=j+1}^n \varphi_2(x,y_i)
\Biggr) \Biggr).
\end{eqnarray*}

Therefore it follows by Theorem~\ref{thmmain2}, $\varphi_i\le1$ for
$i=1,2$ and Campbell's formula [see \citet{dvj08}, Section~9.5] that
\begin{eqnarray*}
d_{\mathrm{TV}}\bigl(\mathscr{L}(\Xi),\mathscr{L}(\mathrm
{H})\bigr) &\leq&
c_1(\lambda) \mathbb{E} \biggl( \int_{\mathcal{X}} \int
_{\mathcal{X}} \beta(x) \bigl| \varphi_1(x,y) -
\varphi_2(x,y) \bigr| \bolds{\alpha}(d x) \Xi(dy) \biggr)
\\
&=& c_1(\lambda) \int_{\mathcal{X}} \int
_{\mathcal{X}} \beta(x) \nu(y) \bigl| \varphi_1(x,y) -
\varphi_2(x,y) \bigr| \bolds{\alpha}(d x) \bolds{\alpha}(dy).
\end{eqnarray*}

Regarding the bound on $c_1(\lambda)$ we have for all $x,y\in\mathcal
{X}$ and $\xi\in\mathfrak{N}$ that $\lambda(x\mvert\xi+\delta
_y)=\lambda(x\mvert\xi)\varphi_2(x,y)$ and $\lambda(x\mvert\xi
)\le\beta(x)$, with equality if $\xi=\varnothing$. Thus by
Theorem~\ref{thmeps-rec},
\begin{eqnarray*}
\varepsilon&=&\sup_{\|\xi-\eta\|=1}\int_\mathcal{X}\bigl|
\lambda(x\mvert\xi)-\lambda(x\mvert\eta)\bigr| \bolds{\alpha}(dx)
\\
&=&\sup_{y\in\mathcal{X}, \xi\in\mathfrak{N}}\int_\mathcal
{X}\bigl|\lambda(x
\mvert\xi+\delta_y)-\lambda(x\mvert\xi)\bigr| \bolds{\alpha} (dx)
\\
&=&\sup_{y\in\mathcal{X}, \xi\in\mathfrak{N}}\int_\mathcal
{X}\lambda(x\mvert
\xi)\bigl|\varphi_2(x,y)-1\bigr| \bolds{\alpha}(dx)
\\
&=& \sup_{y\in\mathcal{X}}\int_\mathcal{X}\beta(x)
\bigl(1-\varphi_2(x,y)\bigr) \bolds{\alpha}(dx).
\end{eqnarray*}
Furthermore, $|\lambda(x\mvert\xi)-\lambda(x\mvert\eta)| \le
\beta(x)$ for all $x\in\mathcal{X}$ and for all $\xi,\eta\in
\mathfrak{N}$. Thus $c \leq\int_\mathcal{X}\beta(x) \bolds
{\alpha}(dx)$.
\end{pf*}

\begin{appendix}
\section*{Appendix: The case of a nondiffuse reference measure~\texorpdfstring{$\alpha$}{alpha}}\label{app}
In the main part we restrict ourselves to a diffuse reference measure
$\bolds{\alpha}$; see Section~\ref{ssecdiffuse}. This appendix shows
that our results remain true for general $\bolds{\alpha}$ (always finite).

Suppose that $\bolds{\alpha}$ is not diffuse. Consider then instead of
$\mathcal{X}$ the extended space $\widetilde{\mathcal{X}}= \mathcal
{X}\times[0,1]$
and equip it with the pseudometric $\tilde{d}((x,u),(y,u)) = d(x,y)$.
In the main part the metric $d$ serves the double purpose of inducing
the $\sigma$-algebra on $\mathcal{X}$ and allowing us to define balls
when it comes to more detailed considerations of conditional
intensities. For the first purpose, which does not require an explicit
metric, we just use the product topology on $\widetilde{\mathcal
{X}}$; for the second
purpose we use the pseudometric $\tilde{d}$. Note that the topology
induced by $\tilde{d}$ is coarser than the product topology, so there
are no measurability problems. Define $\tilde{\bolds{\alpha}}=
\bolds{\alpha}\otimes\Leb
\vert_{[0,1]}$, and denote by $\widetilde{\Po}_1$ the distribution of the
Poisson process with intensity measure $\tilde{\bolds{\alpha}}$.
Note that $\tilde{\bolds{\alpha}}
$ is always a diffuse measure, and that a $\widetilde{\Po}_1$-process is
simple, that is, almost surely free of multi-points.

Transform a point process $\Xi= \sum_{i=1}^N \delta_{X_i}$ on
$\mathcal{X}$ into a point process $\tXXi$ on~$\widetilde{\mathcal
{X}}$ by randomizing
its points in the new coordinate. More precisely let $\tXXi= \sum
_{i=1}^N \delta_{(X_i,U_i)}$, where $U_1, U_2, \ldots$ are i.i.d.
$\Leb
\vert_{[0,1]}$-distributed random variables that are independent
of~$\Xi$. Following \citet{kallenberg86}, we refer to $\tXXi$ as the
\emph{uniform randomization} of $\Xi$. Writing $\widetilde{\mathfrak
{N}}$ for the space
of finite counting measures on~$\widetilde{\mathcal{X}}$, we
introduce the projection
$\pi_{\mathcal{X}} \colon\widetilde{\mathfrak{N}}\to\mathfrak
{N}$, $\pi_{\mathcal
{X}} ( \sum_{i=1}^n \delta_{(x_i,u_i)} ) = \sum_{i=1}^n
\delta_{x_i}$. Note that the image measure of $\widetilde{\Po}_1$ under
$\pi_{\mathcal{X}}$ is $\widetilde{\Po}_1 \pi
_{\mathcal
{X}}^{-1} = \Po_1$. If $\Xi$ is a Gibbs process with density $u$
with respect to $\Po_1$, then a short calculation shows that
$\tXXi$~is a Gibbs process with density
\[
u_{\tXXi}(\tilde{\xi}) = u\bigl(\pi_{\mathcal{X}}(\tilde{\xi
})\bigr)
\qquad\mbox{for any } \tilde{\xi}\in\widetilde{\mathfrak{N}}
\]
with respect to $\widetilde{\Po}_1$.

Hence the conditional intensity of $\tXXi$ is given by
\[
\tilde{\lambda}\bigl((x,u) \mvert\tilde{\xi}\bigr) = \frac
{u(\pi_{\mathcal
{X}}(\tilde{\xi}
)+\delta_x)}{u(\pi_{\mathcal{X}}(\tilde{\xi}))} =
\lambda\bigl(x \mvert\pi_{\mathcal{X}}(\tilde{\xi})\bigr).
\]

Let $\widetilde{\mathcal{F}}_{\mathrm{TV}}$ be the class of
measurable functions $f \colon\widetilde{\mathfrak{N}}\to
[0,1]$ for which $f(\tilde{\xi}) = f(\tilde{\eta})$ whenever\vspace*{1pt} $\pi
_{\mathcal{X}}(\tilde{\xi}) = \pi_{\mathcal{X}}(\tilde{\eta})$.
Furthermore, let $\widetilde{\mathcal{A}}$ be the generator of an SBDP
on $\widetilde{\mathfrak{N}}$ with
birth rate $\tilde{\lambda}(\cdot\mvert\cdot)$ and unit per-capita death
rate. Note that
\[
\tilde{c}_1(\tilde{\lambda})=\sup_{f\in\widetilde{\mathcal
{F}}_{\mathrm{TV}}} \sup
_{\tilde{x}\in\widetilde{\mathcal{X}},
\tilde{\xi}\in
\widetilde{\mathfrak{N}}}\bigl|h_f(\tilde{\xi}+\delta_{\tilde
{x}})-h_f(
\tilde{\xi})\bigr| = c_1(\lambda).
\]

Then for any two-point processes $\Xi$ and $\mathrm{H}$ on $\mathcal{X}$
we can show by slightly adapting the proof of Theorem~\ref{thmmain2} that
\begin{eqnarray*}
d_{\mathrm{TV}} \bigl( \mathscr{L}(\Xi),\mathscr{L}(\mathrm{H})
\bigr) &=& \sup
_{f
\in\widetilde{\mathcal{F}}_{\mathrm{TV}}} \bigl| \mathbb{E}f(\tXXi) -
\mathbb{E}f(\widetilde{
\mathrm{H}}) \bigr|
\\
&=& \sup_{f \in\widetilde{\mathcal{F}}_{\mathrm{TV}}} \bigl| \mathbb
{E}\widetilde{\mathcal{A}}h_f(
\tXXi) \bigr|
\\
&\leq& \tilde{c}_1(\tilde{\lambda}) \int_{\widetilde{\mathcal{X}}}
\mathbb{E} \bigl| \tilde{\nu}\bigl((x,u) \mvert\tXXi\bigr) - \tilde
{\lambda}
\bigl((x,u) \mvert\tXXi\bigr) \bigr| \tilde{\bolds{\alpha}} \bigl(
d(x,u) \bigr)
\\
&=& \tilde{c}_1(\tilde{\lambda}) \int_{\mathcal{X}} \int
_0^1 \mathbb{E} \bigl| \nu\bigl(x \mvert
\pi_{\mathcal{X}}(\tXXi)\bigr) - \lambda\bigl(x \mvert\pi
_{\mathcal{X}}(
\tXXi)\bigr) \bigr| \,du \bolds{\alpha}(dx)
\\
&=& c_1(\lambda) \int_{\mathcal{X}} \mathbb{E} \bigl| \nu(x
\mvert\Xi) - \lambda(x \mvert\Xi) \bigr| \bolds{\alpha}(dx),
\end{eqnarray*}
where for the last equality the two expectations are the same by the
transformation theorem.

The upper bound for $c_1(\lambda)=\tilde{c}_1(\tilde{\lambda})$ in
inequality~(\ref{eqc1}) can be obtained analogously as before by
bounding the expected coupling time between two SBDPs with generator
$\widetilde{\mathcal{A}}$ whose starting configuration differs in
only one point. The
expected coupling time will not be larger than before, because we may
match the additional components $u \in[0,1]$ of any new born points
perfectly in the two processes.

For the more general statement in Remark~\ref{remubermain2} to hold,
we have to replace condition~$(\Sigma)$ by condition $(\Sigma')$; see
\citet{kallenberg86}, Section~13.2.

For the other results in Section~\ref{secdist}, which are essentially
corollaries of Theorem~\ref{thmmain2}, it is easy to verify that we
did not use the fact that $\alpha$ is diffuse, in particular not that
$\xi, \eta\in\mathfrak{N}$ are multi-point free, except for
notational purposes.
\end{appendix}

\section*{Acknowledgements}
We thank the referees for their pertinent comments, which have led to
an improvement of this article.



%

\printaddresses

\end{document}